\documentclass{article} 
\usepackage{texdraw}\usepackage{graphicx}
\usepackage{amsmath} 
\usepackage{color}
\usepackage{colordvi}
\usepackage{texdraw}\usepackage{graphicx}
\usepackage{dcpic,pictex}

\definecolor{purple}{cmyk}{0.45,0.86,0,0}\definecolor{brickred}{cmyk}{0,0.89,0.94,0.28}
\definecolor{maroon}{cmyk}{0,0.87,0.68,0.32}\definecolor{lyellow}{cmyk}{0,0,0.68,0}
\definecolor{magenta}{cmyk}{0,1,0,0}
\definecolor{zals}{cmyk}{0.75,0.12,0.98,0.5}
\definecolor{pamats}{cmyk}{0.54,0.00,0.84,0.0}
\definecolor{brickred}{cmyk}{0,0.89,0.94,0.28}
\definecolor{cadetblue}{cmyk}{0.62,0.57,0.23,0}
\definecolor{teo}{cmyk}{0.63,0.83,0.33,0.13}

\def\newpic#1{%
   \def\emline##1##2##3##4##5##6{%
      \put(##1,##2){\special{em:point #1##3}}%
      \put(##4,##5){\special{em:point #1##6}}%
      \special{em:line #1##3,#1##6}}}
\newpic{}

\usepackage[cp1251,cp1257]{inputenc}\usepackage[english,latvian]{babel}\usepackage{amssymb,latexsym,amsfonts}

\newtheorem{teor}{Theorem}[section]
\newtheorem{theorem}[teor]{Theorem}
\newtheorem{corollary}[teor]{Corollary}
\newtheorem{lemma}[teor]{Lemma}
\newtheorem{proposition}[teor]{Proposition}
\newtheorem{definition}[teor]{Definition}

\newtheorem{example}[teor]{Example}

\input cyracc.def

\begin{document}
\selectlanguage{english}
\title{The Lamplighter Group}
\date{}
\author{J\=anis Buls \\
{\small Department  of  Mathematics, University of Latvia, Jelgavas iela 3,}\\
{\small R\=\i ga, LV-1004 Latvia,
buls@edu.lu.lv
}}
\def\keywords{\begin{center}{\bf Keywords}\end{center}
{automata (machines) groups, the lamplighter group}}

\maketitle
\title{}
\abstract{\small{We are following \cite{ski}. Nevertheless we are interested only in claryfication that the lamplighter group can be realized as a 2--states Mealy machine. }}


\section{ Preliminaries} 

We use standard conventions.
Nethertheless we recall some notation. For more details see \cite{buls}.

Let $A$ be a finite non-empty set and $A^*$ be
the free monoid generated by $A$. The set $A$ is also called an {\em 
alphabet}, its  elements are called {\em letters} and 
those of $A^*$ are called {\em finite words}. The identity element of $A^*$ is called 
an {\em empty word} and denoted by $\lambda$. We set 
$A^+=A^*\backslash\{\lambda\}$.

The length of a word $w$, denoted by $|w|$, is the number of occurrences of the letters
in $w$. In other words, if $w=w_0w_1\cdots w_{n-1}$ with $w_i\in A$, $0\le i<n$, then $|w|=n$. In
particular, the length of the empty word is zero. The set of words of length $k$ over $A$ is denoted by $A^k$.

\begin{definition}
A 2--sorted algebra $\mathcal{M}=\langle Q,A,\circ ,\ast \rangle$ is called  a 
{\em Mealy machine} if $Q,A$ are finite, nonempty sets, the mappings 
$Q\times A \stackrel{\circ}{\longrightarrow}Q$,
$Q\times A \stackrel{\ast }{\longrightarrow}A$
are  total  functions. 
\end{definition}

\begin{itemize}
	\item $Q$ is called  a set of states;
	\item $A$  --- an alphabet;
	\item $\circ$ --- the transition function;
	\item $*$ --- the output function.
\end{itemize}
 One draws 
\begin{eqnarray}
\bigcirc &\xrightarrow{a/b}& \bigcirc \label{for1}\\
q_1 & & q_2 \nonumber
\end{eqnarray}
\medskip\\
to mean $q_1 \circ a= q_2, q_1 *a = b$. That is, in the state $q_1$ with the input $a$, the machine
outputs $b$ and goes to the state $q_2$.

The machine is said to be {\em invertible} if for each state $q\in Q$ the output function is a permutation.
The inverse machine is the
machine $\mathcal{M}^{-1}$ obtained by switching the input $a$ and the output $b$ on each
arrow (\ref{for1}). The resulting transition in the inverse machine is
(\ref{for2}). 

\begin{eqnarray}
\bigcirc &\xrightarrow{b/a}& \bigcirc \label{for2}\\
q_1 & & q_2 \nonumber
\end{eqnarray}

The mappings $\circ$ and $\ast$ may be extended to $Q\times A^*$  by defining
\[
\begin{array}{lr}
q\circ \lambda =q,\quad & q\circ (ua)=(q\circ u)\circ a, \\
q\ast \lambda  =\lambda ,\quad & q\ast (ua)=(q\ast u)\#((q\circ u)\ast a)\,,
\end{array}
\] 
for each $q\in Q$, $(u,a)\in A^*\times A$. Here $\#$ means the concatenation of words. In other words, if 
$u=u_0u_1\cdots u_n$, $v=v_0v_1\cdots v_k$, then $u\#v=uv=u_0u_1\cdots u_nv_0v_1\cdots v_k$.

\section{An Ultrametric}

An (indexed) infinite word $x$ on the alphabet $A$ is any total mapping 
$x\,:\,\mathbb{N}\rightarrow A$. We shall set for any $i\ge0$, $x_i=x(i)$
and write
\[
x=(x_i)=x_0x_1\cdots x_n\cdots \;.
\]
The set of all the infinite words over $A$ is denoted by $A^\omega$. 

Let $A^\infty=A^*\cup A^\omega$. A word $u\in A^*$ is a {\em  prefix} of 
$x\in A^\infty$ if there exist $u\in A^*$, $y\in A^\infty$ such that $x=uy$.
The word $y$ is called the {\em suffix} in this situation.
The set  of
 prefixes of $x$ is denoted by $\mathcal{P}(x)$.

\begin{definition} 
A mapping $d:A^\infty \times A^\infty \to \mathbb{R}$ is called the prefix metric, if
\[
d(x,y)=\begin{cases} 0,& {\rm if}\; x=y;\\
2^{-n}, & {\rm if}\; x\ne y \; {\rm and}\; n=\max\{|u|\,|\,u\in \mathcal{P}(x)\cap\mathcal{P}(y)\}.
\end{cases}
\]
\end{definition}

 An  {\em ultrametric} is a metric which satisfies the following strengthened version of the triangle inequality:
\[
d(x,z)\le\max(d(x,y),d(y,z)) \quad {\rm for \; all} \quad x,y,z.
\]

\begin{proposition}
The prefix metric is an ultrametric.
\end{proposition}

$\Box$ (i) If $x=z$, then $d(x,z)=0\le \max(d(x,y),d(y,z))$.

(ii) If $x\ne z$, then $d(x,z)=2^{-n}$, where $\forall i<n\;(x_i=z_i)$ but $x_n\ne z_n$.
Hence, $x_n\ne y_n$ or $y_n\ne z_n$. Therefore $d(x,y)\ge 2^{-n}$ or $d(y,z)\ge 2^{-n}$. Thus,
$d(x,z)=2^{-n}\le \max(d(x,y),d(y,z))$.
\rule{2mm}{2mm}

 If $\mathcal{M}=\langle Q,A,\circ ,\ast \rangle$ is   a 
 Mealy machine, then the mapping $\ast$ may be extended to $Q\times A^\infty$  by defining
\[
q*x=\lim\limits_{n\to \infty}q*x[0,n],
\]
where $x[0,n]=x_0x_1\cdots x_n$.

We like to be more flexible therefore we sometimes write $uq$ instead of $q*u$ and state $q$ identify with map $q*u$. In other words we identify map $\bar q: A^* \to A^*: u\mapsto q*u$ with $q$. 
For an invertible Mealy machine $\mathcal{M}$, the group generated by $\{ \bar q\,|\, q\in Q\}$
 under the operation of composition is called the machine (automaton) group
$\Gamma(\mathcal{M})$.

\section{Sequential functions}

Now we recall standard results about sequential functions. Proofs see in \cite{buls}.

\begin{definition} A total mapping $f : A^* \to B^*$ is called a sequential function if

{\rm (i)} $\forall u\in A^* \; |u| = |f(u)|$;

{\rm (ii)} $u \in \mathcal{P}(v) \Rightarrow f(u) \in \mathcal{P}(f(v))$.
\end{definition}

\begin{corollary} For all sequential functions, we have that if
\[
u \in \mathcal{P}(v) \cap \mathcal{P}(w),
\]
then
\[
f(u) \in \mathcal{P}(f(v)) \cap \mathcal{P}(f(w)).
\]
\end{corollary}

\begin{definition} Let $f : A^* \to B^* $ be a sequential function and \linebreak
$u \in A^*$, then $f_u(v)$ defines a suffix of mapping $f(uv)$ with length $|v|$.
The mapping $f_u$ is called a remainder of sequential function $f$.
\end{definition}

\begin{corollary} $f(uv) = f(u)f_u(v)$.
\end{corollary}

\begin{lemma} $uw \in \mathcal{P}(uv) \Rightarrow w \in \mathcal{P}(v)$.
\end{lemma}

\begin{proposition} \label{P3.6} The remainder $f_u$ is a sequential function.
\end{proposition}

\begin{lemma} If mapping $f : A^* \to B^*$ is a sequential function,
then
\[
f_u(vw) = f_u(v)f_{uv}(w).
\]
\end{lemma}

\begin{proposition} Let a mapping $f : A^* \to B^*$ be a sequential
function. If $f_u = f_{u'}$, then
\[
\forall v\in A^* \; f_{uv} = f_{u'v}.
\]
\end{proposition}

\begin{lemma}
If $f : A^* \to B^*$ is a a sequential
function, then \[f_{uv}=(f_u)_v.\]
\end{lemma}

$\Box$ Let $w\in A^*$, then
\begin{eqnarray*}
f(uvw) &=& f(uv)f_{uv}(w)= f(u)f_u(v)f_{uv}(w),\\
f(uvw) &=& f(u)f_u(vw)=f(u)f_u(v)(f_u)_v(w).
\end{eqnarray*}
Thus, $f_{uv}(w)=(f_u)_v(w)$.
\rule{2mm}{2mm}

\begin{proposition} Let mappings $f : A^* \to B^*$ and $g : B^* \to C^*$
be sequential functions, then
\[
\forall u\in A^* \; (g \circ f)_u = g_{f(u)} \circ f_u.
\]
Here $(g\circ f)(x)=f(g(x))$.
\end{proposition}

\begin{definition} Let a mapping $f : A^* \to B^*$ be a sequential function. The function $f$ defines the set
\[
Q_f = \{f_u\,|\,u \in A^*\},
\]
where $f_u$ is a remainder of $f$. The function $f$ is called a restricted
sequential function (or a sequential function on a restricted domain) if
the set $Q_f$ is finite.
\end{definition}

\begin{definition}
A 2-sorted algebra $\mathcal{M}=\langle Q,A, q_0, \circ ,\ast \rangle$ is called  an 
{\em initial Mealy machine} if 
\begin{itemize}
	\item $q_0\in Q$;
	\item $\langle Q,A, \circ ,\ast \rangle$ is a Mealy machine.
\end{itemize}
\end{definition}

\begin{theorem} \label{teor3.13}
A function $f : A^* \to A^*$ is the restricted
sequential function if and only if there exists an
initial Mealy machine $\mathcal{M} = \langle Q , A,  q_0, \circ, \ast\rangle$ such that
\[
\forall v \in A^* \quad f(v) = q_0 * v.
\]
\end{theorem}

\section{$\omega$--sequential functions}

\begin{definition} A total mapping $\zeta : A^\omega \to B^\omega$ is called an \linebreak
$\omega$--sequen\-tial function if
it satisfies such condition:
\[
u \in \mathcal{P}(x)\cap \mathcal{P}(y) \Rightarrow \zeta(x)[0,|u|) = \zeta(y)[0,|u|).
\]
\end{definition}

Here, if $\zeta(x)=z_0z_1\cdots z_n\cdots$ then $\zeta(x)[0,n)=z_0z_1\cdots z_{n-1}$. 

\begin{proposition}
If $\zeta: A^\omega \to B^\omega$, $\xi: B^\omega \to C^\omega$ are  $\omega$--sequential func\-tions then
$\zeta\xi=\xi\circ \zeta: A^\omega \to C^\omega$ is an $\omega$--sequential function.
\end{proposition}

$\Box$ Let $u\in \mathcal{P}(x)\cap \mathcal{P}(y)$ then $\zeta(x)[0,|u|) = \zeta(y)[0,|u|)$.
Hence \[\xi(\zeta(x))[0,|u|)=\xi(\zeta(y))[0,|u|). \quad
\rule{2mm}{2mm}\]

If $\zeta: A^\omega \to B^\omega$ is an $\omega$--sequential function then 
\[
\forall x\in A^\omega\; \forall y\in A^\omega \; \forall u\in A^*\quad \zeta(ux)[0,|u|)=\zeta(uy)[0,|u|).
\]
Hence $\zeta(ux)=\epsilon(u)\# \zeta_u(x)$, where $\epsilon$ and $\zeta_u$ are mappings
\[
\epsilon: A^{|u|} \to B^{|u|}, \quad \zeta_u: A^\omega \to B^\omega.
\]
The mapping $\zeta_u$ is called the {\em remainder} of $\zeta$ for the word $u$.

\begin{proposition}
$\zeta_u$ is an $\omega$--sequential function.
\end{proposition}

$\Box$ Let $v\in \mathcal{P}(x)\cap \mathcal{P}(y)$ then $uv\in \mathcal{P}(ux)\cap \mathcal{P}(uy)$.
Since $\zeta$ is an $\omega$--sequential function then
\begin{eqnarray*}
\epsilon(u)\#\zeta_u(vx)[0,|v|) &=&\zeta(uvx)[0,|uv|)=\zeta(uvy)[0,|uv|)\\
&=&\epsilon(u)\# \zeta_u(vy)[0,|v|) \quad \rule{2mm}{2mm}
\end{eqnarray*}

\begin{lemma}
If $\zeta: A^\omega \to B^\omega$ is an $\omega$--sequential function then \\
$\zeta_{uv}=(\zeta_u)_v$.
\end{lemma}

\begin{eqnarray*}
\Box \quad \zeta(uvx) &=& \epsilon_1(u)\#\zeta_u(vx), \quad {\rm where}\\
\epsilon_1(u) &=&  \zeta(uvx)[0,|u|);\\
\zeta_u(vx) &=& \epsilon_2(v)\# (\zeta_u)_v(x), \quad {\rm where}\\
\epsilon_2(v) &=& \zeta_u(vx)[0,|v|);\\
\zeta(uvx) &=& \epsilon_3(uv)\# \zeta_{uv}(x), \quad {\rm where} \\
\epsilon_3(uv) &=& \zeta(uvx)[0.|uv|).
\end{eqnarray*}
Since $|\epsilon_1(u)|+|\epsilon_2(v)|=|u|+|v|=|uv|=\epsilon_3(uv)$ then
$(\zeta_u)_v=\zeta_{uv}$. \rule{2mm}{2mm}

\begin{lemma}
If $\zeta: A^\omega \to A^\omega$ is a bijective $\omega$--sequential function then 
$\epsilon: A^{|u|} \to A^{|u|}$ is a bijection.
\end{lemma}

$\Box$   Let $u\ne v$. If $\epsilon(u)=\epsilon(v)$ then there exists $w'\in A^{|u|}$ such that
\[
\forall w\in A^{|u|}\quad \epsilon(w)\ne w'
\]
because the set $A^{|u|}$ is finite.

Let $w'y\in A^\omega$ then there exists  $x\in A^{\omega}$  such that $\zeta(x)=w'y$ because $\zeta$ is a surjection.
Let $w\in\mathcal{P}(x)$ and $|w|=|w'|$ then $x=wx'$ for some $x'\in A^\omega$ and
\[
w'y=\zeta(x)=\zeta(wx')=\epsilon(w)\#\zeta_w(x').
\]
Therefore $\epsilon(w)=w'$. Contradiction! 

This means that $\epsilon: A^{|u|} \to A^{|u|}$ is an injection.
Now we can conclude that $\epsilon: A^{|u|} \to A^{|u|}$ is a bijection because the set $A^{|u|}$ is finite.
\rule{2mm}{2mm}

\begin{lemma}
If $\zeta: A^\omega \to A^\omega$ is a bijective $\omega$--sequential function then $\zeta_u$ is a bijection.
\end{lemma}

$\Box$ We know $\zeta(ux)=\epsilon(u)\#\zeta_u(x)$. Let $\epsilon(u)=v$. Since $\epsilon: A^{|u|} \to A^{|u|}$ is a bijection then there exists only one word $u'$ for which $\epsilon: u'\mapsto v$, namely $u$ is such the word.
Therefore the words $vy$ can be obtained only from words $ux$, namely, 
\[
\forall y\in A^\omega\; \exists x\in A^\omega \quad vy=\zeta(ux)=\epsilon(u)\#\zeta_u(x)=v\zeta_u(x)
\]
Hence $y=\zeta_u(x)$. This means that $\zeta_u$ is a surjection.

Since $\zeta$ is an injection then
\[
\epsilon(u)\#\zeta_u(x)=\zeta(ux)\ne \zeta(ux')=\epsilon(u)\#\zeta_u(x')
\]
if $x\ne x'$. Hence $\zeta_u(x)\ne \zeta_u(x')$. This means that $\zeta_u$ is an injection.

Now we can conclude that $\zeta_u$ is a bijection because $\zeta_u$ is both a surjection and an injection.
\rule{2mm}{2mm}

\begin{lemma}
If $\zeta: A^\omega \to A^\omega$ is a bijective $\omega$--sequential function then $\zeta^{-1}$ is a bijective $\omega$--sequential function.
\end{lemma}

$\Box$ Since $\zeta$ is a bijection then $\zeta^{-1}$ exists and it is a bijection. We must prove that
$\zeta^{-1}(wx)[0,|w|)=\zeta^{-1}(wy)[0,|w|)$ for all $w\in A^*$ and all infinite words $x,y$ over alphabet $A$.

Let $ux'=\zeta^{-1}(wx)$, $vy'=\zeta^{-1}(wy)$ and $|u|=|v|=|w|$
then 
\[
\zeta(ux')=\zeta\zeta^{-1}(wx)=wx,\quad  \zeta(vy')=\zeta\zeta^{-1}(wy)=wy
\] and 
$\zeta^{-1}(wx)[0,|w|)=u$, $\zeta^{-1}(wy)[0,|w|)=v$. 
\[
wx=\zeta(ux')=\epsilon(u)\#\zeta_u(x'), \quad wy=\zeta(vy')=\epsilon(v)\#\zeta_v(y').
\]
Hence $\epsilon(u)=w=\epsilon(v)$ and thus $u=v$ because $\epsilon:A^{|u|}\to A^{|u|}$ is a bijection.
Therefore $\zeta^{-1}(wx)[0,|w|)=u=v=\zeta^{-1}(wy)[0,|w|)$.
\rule{2mm}{2mm}

\section{Endomorphisms}

\begin{example} Let's choose an alphabet $A = \{a_1, a_2, \ldots , a_p\}$.
Let's define a p–regular rooted tree $\mathcal{T}_A$ as:
\begin{itemize}
	\item the words of $A^*$ are vertices of the tree;
	\item the empty word $\lambda$ is a root of the tree;
	\item the set $\mathcal{E}_A=\{(u, ua) \,|\, u\in A^* \wedge  a \in A\}$ is the set of arcs.
\end{itemize}
\end{example}

We shall restrict our attention only on such trees $\mathcal{T}_A$.
Illustration when $A=\{0,1\}$ see in the figure \ref{lamp1}

\begin{figure}[h]
{\input{aut47.pic}}
\caption{A rooted tree.}
\label{lamp1}
\end{figure}

\begin{definition} A mapping $f : A^* \to A^*$ is called an endomorphism of a tree $\mathcal{T}_A$ if
\begin{itemize}
	\item $f(\lambda)=\lambda$;
	\item  $ (f(u), f(v)) \in \mathcal{E}_A$  for each arc
$(u, v)\in\mathcal{E}_A$.
\end{itemize}
\end{definition}

If $f$ is a bijection then an endomorphism $f$ is called an {\em automorphism} of the tree $\mathcal{T}_A$.
 We denote by ${\rm End}(\mathcal{T}_A)$ the set of all endomorphisms of the tree $\mathcal{T}_A$ and by 
${\rm Aut}(\mathcal{T}_A)$ we denote the set of all automorphisms of the tree $\mathcal{T}_A$.

\begin{proposition}\label{P5.3}
 A mapping $f : A^* \to A^*$ is a sequential function, if and only if $f \in  {\rm End}(\mathcal{T}_A).$
\end{proposition}

$\Box$ The proof of this fact can be found for example in \cite{buls}. \rule{2mm}{2mm}

\begin{proposition}
If $f\in {\rm Aut}(\mathcal{T}_A)$ then $\forall u\in A^* \; f_u\in{\rm Aut}(\mathcal{T}_A)$.
\end{proposition}

$\Box$ We know (\ref{P3.6}. Proposition and \ref{P5.3}. Proposition ) that \\
$f_u\in {\rm End}(\mathcal{T}_A) $. So we must prove that $f_u$ is a bijection.

Let $f(u)=v$.

(i) Let $w\in A^*$ then exists $w'$ such that $f(w')=vw$. Choose $u'\in\mathcal{P}(w')$ such that $|u'|=|u|$ then $w'=u'w''$ for some $w''$. Hence
\[
vw=f(w')=f(u'w'')=f(u')f_{u'}(w'').
\]
Since $f$ is a bijection then $u'=u$. Therefore $vw=f(u)f_u(w'')$. This means that $f_u$ is a surjection.

(ii) Let $w\ne w'$ then $f(uw)\ne f(uw')$. Hence
\[
f(u)f_u(w)=f(uw)\ne f(uw')=f(u)f_u(w').
\]
This demonstrates that $f_u(w)\ne f_u(w')$. Therefore $f_u$ is an injection.
\rule{2mm}{2mm}
\medskip

If $v\in A^+$ then we denote by $v^\omega$ the infinite word
\[
v^\omega=vv\cdots v\cdots 
\]

\begin{proposition}\label{prop5.5}
Every $\omega$--sequential function $\zeta: A^\omega \to B^\omega$ induces the sequential function
$\breve\zeta: A^* \to B^*$.
\end{proposition}

$\Box$ At first we must clarify:

--- What does it mean induces?
\begin{itemize}
	\item $\breve\zeta(\lambda)=\lambda$;
	\item $\breve\zeta(u)=\zeta(ua^\omega)[0,|u|)$ for all $(u,a)\in A^*\times A$. 
\end{itemize}

Now we shall prove that $\breve\zeta$ is a sequential function. By the definition of $\breve\zeta$ we can conclude
\[
\forall u\in A^*\quad |u|=|\breve\zeta(u)|.
\]

Let $u\in\mathcal{P}(v)$ then
\begin{eqnarray*}
\breve\zeta(u) &=& \zeta(ua^\omega)[0,|u|)\\
&=& \zeta(va^\omega)[0,|u|)\in \mathcal{P}(\zeta(va^\omega))[0,|v|)=\mathcal{P}(\breve\zeta(v)). \quad
\rule{2mm}{2mm}
\end{eqnarray*}

\begin{proposition}\label{P5.5}
If an $\omega$--sequential function $\zeta: A^\omega \to A^\omega$ is a bijection then the induced sequential function
$\breve\zeta\in {\rm Aut}(\mathcal{T}_A)$.
\end{proposition}

$\Box$ Let $v\in A^*$. Since $\zeta$ is the bijection then there exists $ux\in A^\omega$ such that $|u|=|v|$ and
$\zeta(ux)=va^\omega$ for some $a\in A$. Hence
\[
v=\zeta(ux)[0,|u|)=\zeta(ua^\omega)[0,|u|)=\breve\zeta(u).
\]

Therefore $\forall n\; \breve\zeta: A^n \to A^n$ is a surjection. Since alphabet $A$ is a finite set then 
$\forall n\; \breve\zeta: A^n \to A^n$ is a bijection. This means that $ \breve\zeta: A^* \to A^*$ is a bijection.
\rule{2mm}{2mm}

\begin{corollary}
If an $\omega$--sequential function $\zeta: A^\omega \to A^\omega$ is a bijection then the induced sequential function
$\breve\zeta\in {\rm Aut}(\mathcal{T}_A)$.
\end{corollary}

$\Box$ \ref{P5.5}. Proposition, ${\rm Aut}(\mathcal{T}_A)$ definition and  \ref{P5.3}. Proposition.
\rule{2mm}{2mm}

\begin{proposition}
If $f$ is the induced function of an $\omega$--sequential function $\zeta$ and $\breve f_u$ is the induced function of the $\omega$--sequential function $\zeta_u$ then $f_u=\breve f_u$.
\end{proposition}

$\Box$ Let $\zeta: A^\omega \to B^\omega$ be an $\omega$--sequential function and $v\in A^*$.
 By the definition of $f$ we can conclude:
\begin{eqnarray*}
f(uv) &=& f(u)\#f_u(v),\\
f(uv) &=& \zeta(uvx)[0,|uv|), \\
\zeta(uvx) &=& \epsilon(u)\# \zeta_u(vx), \quad |\epsilon(u)|=|u|, \\
\breve f_u(v) &=& \zeta_u(vx)[0,|v|), \\
\zeta(uvx)[0,|uv|) &=&  \epsilon(u)\# \zeta_u(vx)[0,|v|)=\epsilon(u)\#\breve f_u(v).
\end{eqnarray*}

Hence
\[
f(u)\# f_u(v)=\zeta(uvx)[0,|uv|)=\epsilon(u)\#\breve f_u(v).
\]
Since $|\epsilon(u)|=|u|=|f(u)|$ then $f_u(v)=\breve f_u(v)$.
\rule{2mm}{2mm}

\begin{proposition}\label{prop5.8}
Every sequential function $f: A^* \to B^*$ indu\-ces the $\omega$--sequential function
$\hat f: A^\omega \to B^\omega$.
\end{proposition}

$\Box$ At first we must clarify:

--- What does it mean induces?
\[
\hat f(x)=\lim\limits_{n\to \infty} f(x[0,n))
\]

Since $f: A^* \to B^*$ is sequential function then \\
$f(x[0,n))\in\mathcal{P}(f(x[0,n+k)))$ for all natural $k$. Terefore $\lim\limits_{n\to \infty} f(x[0,n))$ exists and $\hat f(x)[0,n)=f(x[0,n))$. Hence
\[
\hat f(ux)[0,|u|)=f(u)=\hat f(uy)[0,|u|)
\]
for all $x,y\in A^\omega$. This means that $\hat f: A^\omega \to B^\omega$ is an $\omega$--sequential function.
\rule{2mm}{2mm}

\section{Formal power series}

Let $R$ be an  associative commutative ring with a multiplicative identity element and $X\not\in R$. We shall identify $R^\omega$ with $R[[X]]$ via
\[
(a_0,a_1,a_2,\ldots, a_n, \ldots) \mapsto \sum\limits_{k=0}^\infty a_kX^k,
\]
 where on the left hand side we write sometimes elements of $R^\omega$ as tuples to avoid
confusion between concatenation of words and the multiplication in $R$.
If $f=\sum\limits_{k=0}^\infty a_kX^k$ then we write $f(n)=a_n$.

\begin{definition}
Algebra $\langle R[[X]],+, \cdot\rangle$ is called a formal power series if
\begin{eqnarray*}
\sum\limits_{k=0}^\infty a_kX^k + \sum\limits_{k=0}^\infty b_kX^k &=& \sum\limits_{k=0}^\infty (a_k+b_k)X^k,\\
\left(\sum\limits_{k=0}^\infty a_kX^k\right)\left(\sum\limits_{k=0}^\infty b_kX^k\right) &=& \sum\limits_{k=0}^\infty \left(\sum\limits_{i=0}^k a_ib_{k-i}\right)X^k
\end{eqnarray*}
\end{definition}

\begin{proposition}
$R[[X]]$ is an  associative commutative ring with a multiplicative identity element $10^\omega$.
\end{proposition}

$\Box$ Let $f=\sum\limits_{k=0}^\infty a_kX^k,$ $g=\sum\limits_{k=0}^\infty b_kX^k,$ $h=\sum\limits_{k=0}^\infty a_kX^k$ then 

\begin{eqnarray*}
fg &=& \sum\limits_{i=0}^\infty \sum\limits_{j=0}^i a_j b_{i-j}X^i,\\
(fg)h &=& \sum\limits_{\varkappa=0}^\infty \sum\limits_{i=0}^\varkappa \sum\limits_{j=0}^i a_j b_{i-j} c_{\varkappa-i}X^\varkappa,\\
gh &=& \sum\limits_{s=0}^\infty \sum\limits_{n=0}^s b_n c_{s-n}=\sum\limits_{s=0}^\infty \beta_sX^s,\quad {\rm where}\\
\beta_s &=& \sum\limits_{n=0}^s b_n c_{s-n},\\
f(gh) &=& \sum\limits_{\varkappa=0}^\infty \sum\limits_{k=0}^\varkappa a_k \beta_{\varkappa-k}X^\varkappa=
\sum\limits_{\varkappa=0}^\infty \sum\limits_{k=0}^\varkappa a_k \sum\limits_{n=0}^{\varkappa-k} b_n c_{\varkappa-k-n}X^\varkappa\\
&=& \sum\limits_{\varkappa=0}^\infty \sum\limits_{k=0}^\varkappa  \sum\limits_{n=0}^{\varkappa-k}a_k b_n c_{\varkappa-k-n}X^\varkappa
\end{eqnarray*}
We must prove that sums
\begin{eqnarray}\label{f3}
\sum\limits_{\varkappa=0}^\infty \sum\limits_{i=0}^\varkappa \sum\limits_{j=0}^i a_j b_{i-j} c_{\varkappa-i}X^\varkappa,
\end{eqnarray}
\begin{eqnarray}\label{f4}
\sum\limits_{\varkappa=0}^\infty \sum\limits_{k=0}^\varkappa  \sum\limits_{n=0}^{\varkappa-k}a_k b_n c_{\varkappa-k-n}X^\varkappa
\end{eqnarray}
are equal. Let's choose $k=j$. Since $j\le i$ then $i-j\ge 0$. Notice $k=j$ and $i\le\varkappa$ therefore
$\varkappa - k=\varkappa-j\ge i-j$. Let $n=i-j$ then
\[
a_kb_nc_{\varkappa-k-n}=a_jb_{i-j}c_{\varkappa-i}
\]
is from sum  (\ref{f4}). Since the number of summands is equel in both sums (\ref{f3}) and (\ref{f4}) then
(illustration see in the table bellow) this means that the sum (\ref{f3}) is equel to the sum (\ref{f4}). 

\begin{center}
\begin{tabular}{l|llll||l|lllll}\label{t1}
$i$&$j$&&&& $k$ & $n$\\
\hline
0&0&&&& 0 & 0& 1& $\ldots$ &$\varkappa-1$ & $\varkappa$ \\ 
1 &0 &1 &&& 1 &0 &1 &$\ldots$ &$\varkappa-1 $\\
$\cdot$ &$\cdot$ &$\cdot$ & $\cdot $ & &$\cdot$ &$\cdot$ &$\cdot$ &$\ldots $  \\
$\cdot$ &$\cdot$ &$\cdot$ &$\cdot \, \cdot $ & &$\cdot$ &$\cdot$ &$\cdot$ &$\cdot$\\
$\cdot$ &$\cdot$ &$\cdot$ &$\ldots$  & &$\cdot$& $\cdot$ &$\cdot$ \\
$\varkappa$ &0 &1 &$\ldots$ &$\varkappa$ &$\varkappa$ &0
\end{tabular}
\end{center}

Thus we have proved associative law, namely $(fg)h=f(gh)$.

\[
\sum\limits_{j=0}^i a_jb_{i-j}=a_0b_i+a_1b_{i-1}+\ldots+ a_{i-1}b_1+a_ib_0=\sum\limits_{j=0}^i b_ja_{i-j}
\]
Hence
\[
fg=\sum\limits_{i=0}^\infty\sum\limits_{j=0}^i a_jb_{i-j}X^i=\sum\limits_{i=0}^\infty\sum\limits_{j=0}^i b_ja_{i-j}X^i=gf.
\]
This is commutative law.
\medskip

$10^\omega\mapsto 1 X^0+\sum\limits_{i=1}^\infty 0X^i$. 
Hence
\begin{eqnarray*}
(1 X^0+\sum\limits_{i=1}^\infty 0X^i)f &=& 
(1 X^0+\sum\limits_{i=1}^\infty 0X^i)\sum\limits_{k=0}^\infty a_kX^k \\
&=& \sum\limits_{k=0}^\infty a_kX^k=
\left(\sum\limits_{k=0}^\infty a_kX^k\right)(1 X^0+\sum\limits_{i=1}^\infty 0X^i)\\
&=&f(1 X^0+\sum\limits_{i=1}^\infty 0X^i)
\end{eqnarray*}
Thus $(1 X^0+\sum\limits_{i=1}^\infty 0X^i)$ is the multiplicative identity element in $R[[X]]$. We idetify $r X^0+\sum\limits_{i=1}^\infty 0X^i$ with $r\in R$. So the multiplicative identity element of the ring $R$ we idetify as  the multiplicative identity element of $R[[X]]$ too.

\begin{eqnarray*}
f(g+h) &=& \left(\sum\limits_{k=0}^\infty a_kX^k\right)\left(\sum\limits_{k=0}^\infty b_kX^k+\sum\limits_{k=0}^\infty c_kX^k\right)\\
&=&\left(\sum\limits_{k=0}^\infty a_kX^k\right)\left(\sum\limits_{k=0}^\infty (b_k+c_k)X^k\right)\\
&=& \sum\limits_{k=0}^\infty \sum\limits_{j=0}^k a_j(b_{k-j}+c_{k-j}) X^k\\
&=& \sum\limits_{k=0}^\infty \sum\limits_{j=0}^k(a_jb_{k-j}+a_jc_{k-j})X^k\\
&=& \sum\limits_{k=0}^\infty \left(\sum\limits_{j=0}^ka_jb_{k-j}+\sum\limits_{j=0}^ka_jc_{k-j}\right)X^k\\
&=&  \sum\limits_{k=0}^\infty \sum\limits_{j=0}^ka_jb_{k-j}X^k+ \sum\limits_{k=0}^\infty \sum\limits_{j=0}^ka_jc_{k-j}X^k\\
&=& fg+fh
\end{eqnarray*}
This is distributive law.

Note 0 is an additive identity and $-f=\sum\limits_{k=0}^\infty (-a_k)X^k$ is the additive inverse of $f$. Now we can conclude from definition of the sum in $R[[X]]$ that $\langle R[[X]], +\rangle$ is a commutative group. 
\rule{2mm}{2mm}

\begin{proposition}
Formal power series $f=\sum\limits_{k=0}^\infty a_kX^k$ is invertible in $R[[x]]$ if and only if then $a_0$ is invertible in $R$.
\end{proposition}

$\Box\; \Rightarrow$ Let $g=\sum\limits_{k=0}^\infty b_kX^k$. If $fg=1$ then $a_0b_0=1$. Thus $a_0$ is invertible in $R$.

$\Leftarrow$  If   $a_0$ is invertible in $R$ then $b_0=a_0^{-1}$. If we like to get $a_0b_1+a_1b_0=0$ then
$b_1=-a_0^{-1}a_1b_0$.

The rest is induction. If we like to get $a_0b_n+a_1b_{n-1}+\ldots +a_nb_0=0$ 
then 
\[
b_n=-a_0^{-1}\left( \sum\limits_{k=1}^n a_k b_{n-k}\right). \quad \rule{2mm}{2mm}
\]

\begin{example}
$X$ is not invertible in $R[[X]]$.
\end{example}

\section{Bijections in $R[[X]]$}

For any power series $f(X)\in R[[X]]$ we define \cite{ski} two mappings of $R[[x]]$ given by
\begin{eqnarray*}
\mu[f] &:& g(X) \mapsto f(X)g(X),\\
\alpha[f] &:& g(X) \mapsto f(X) + g(X).
\end{eqnarray*}

\begin{proposition}
$\alpha[f]$ and $\mu[f]$ are $\omega$--sequential functions.
\end{proposition}

$\Box$ If \quad $f=\sum\limits_{k=0}^\infty d_kX^k$, \quad $g=\sum\limits_{k=0}^{n-1} a_kX^k+\sum\limits_{k=n}^\infty b_kX^k$ \quad and \\
$h=\sum\limits_{k=0}^{n-1} a_kX^k+\sum\limits_{k=n}^\infty c_kX^k$ \quad then
\begin{eqnarray*}
\alpha[f](g)[0,n) &=& \sum\limits_{k=0}^n (d_k+a_k)X^k=\alpha[f](h)[0,n),\\
\mu[f](g)[0,n) &=& \sum\limits_{k=0}^n \left(\sum\limits_{i=0}^k d_i a_{k-i}\right)X^k=\mu[f](h)[0,n). \quad \rule{2mm}{2mm}
\end{eqnarray*}

\begin{proposition}
If $\langle G, \odot\rangle$ is a group then a mapping
\[
T_a:x\mapsto a\odot x
\]
is a bijection.
\end{proposition}

$\Box$ Proof see for example in \cite{We}. \rule{2mm}{2mm}

\begin{corollary}
$\alpha[f]$ is bijective.
\end{corollary}

Let $f_\alpha$ be the sequential function induced by the $\omega$--sequential function $\alpha[f]$ then $f_\alpha\in {\rm Aut}(\mathcal{T}_R)$.

\begin{proposition}
If $\langle G, +, \cdot \rangle$ is a ring and $a\in G$ is a unit then a mapping
\[
T'_a: x\mapsto ax
\]
is a bijection.
\end{proposition}

$\Box$
Let $ax=by$ then $0=ax-ay=a(x-y)$. Since $a$ is unit then $x-y=0$. Thus $x=y$ and so $T'_a$
 is injective. Since $T'_a(a^{-1}y)=aa^{-1}y=y$ then $T'_a$ is 
surjective. 
\rule{2mm}{2mm}

\begin{corollary}
If $f$ is invertible in $R[[x]]$ then  $\mu[f]$ is bijective.
\end{corollary}

Let $f_\mu$ be the sequential function induced by the $\omega$--sequential function $\mu[f]$ then $f_\mu\in {\rm End}(\mathcal{T}_R)$. If $f$ is invertible in $R[[x]]$ then $f_\mu\in {\rm Aut}(\mathcal{T}_R)$.

\begin{proposition}\label{prop7.6}
If $f$ is invertible in $R[[x]]$ then $(\mu[f])^{-1}=\mu[f^{-1}]$.
\end{proposition}

$\Box$ Let $g\in R[[X]]$ then 
\begin{eqnarray*}
\mu[f^{-1}](\mu[f](g)) &=& \mu[f^{-1}](fg)=f^{-1}fg=g,\\
 \mu[f](\mu[f^{-1}](g)) &=& \mu[f](f^{-1}g)=ff^{-1}g=g.
\end{eqnarray*}
This means that $\mu[f^{-1}]=(\mu[f])^{-1}$.
\rule{2mm}{2mm}

\begin{proposition}\label{prop7.7}
If $f$ is invertible in $R[[x]]$ then 
\[\mu[f^{-1}]\alpha[h]\mu[f]=\alpha[fh].
\]
\end{proposition}
\begin{eqnarray*} \Box \quad
g\mu[f^{-1}]\alpha[h]\mu[f] &=& f^{-1}g\alpha[h]\mu[f]= (h+f^{-1}g)\mu[f]=fh+g\\
&=& g\alpha[fh] \quad \rule{2mm}{2mm}
\end{eqnarray*} 

\begin{definition}
For any power series $f(X) = \sum\limits_{k=0}^\infty a_kX^k$ we define also the shift of $f$ by
\[
\sigma(f) =\sum\limits_{k=0}^\infty a_{k+1}X^k
\]
\end{definition}

\begin{corollary}
$f= a_0+\sigma(f)X$.
\end{corollary}

\begin{lemma}
$(1-aX)^{-1}=\sum\limits_{k=0}^\infty a^kX^k$.
\end{lemma}

\begin{eqnarray*}\Box \quad (1-ax)\sum\limits_{k=0}^\infty a^kX^k &=& \sum\limits_{k=0}^\infty a^kX^k-\sum\limits_{k=0}^\infty a^{k+1}X^{k+1}\\
&=& 1+\sum\limits_{k=1}^\infty a^kX^k-\sum\limits_{k=1}^\infty a^kX^k=1 \quad
\rule{2mm}{2mm} \end{eqnarray*}

\begin{lemma}
If $f=\displaystyle\frac{1}{1-aX}$ then $\sigma(f)=af$.
\end{lemma}

\begin{eqnarray*}
\Box \quad \sigma(f) &=& \sum\limits_{k=0}^\infty a^{k+1}X^k,\\
af &=& a\sum\limits_{k=0}^\infty a^kX^k=\sum\limits_{k=0}^\infty a^{k+1}X^k. \quad
\rule{2mm}{2mm}
\end{eqnarray*}

\begin{definition} Let a mapping $\zeta : A^\omega \to B^\omega$ be an $\omega$--sequential function. The function $\zeta$ defines the set
\[
Q_\zeta = \{\zeta_u\,|\,u \in A^*\},
\]
where $\zeta_u$ is a remainder of $\zeta$.  If
the set $Q_f$ is finite then $\zeta$ is called the finite state.
\end{definition}

\begin{lemma}\label{lemma7.13}
If $\zeta\in Q\subseteq Q_\zeta$ and for all $a\in A$ and for every $\xi\in Q$ we have: $\xi_a\in Q$ then $Q=Q_\zeta$.
\end{lemma}

$\Box$ Since $\zeta\in Q$ then for every  $a\in A$ we have $\zeta_a\in Q$.

The rest is induction by length $|w|$ of words $w\in A^*$. Let $|v|=|w|+1$ then $v=wa$ for some $a\in A$.
Since $\zeta_w\in Q$ then $\zeta_v=\zeta_{wa}=(\zeta_w)_a\in Q$.
\rule{2mm}{2mm}

\begin{theorem}\label{teor7.14}
If $\displaystyle f(x)=\frac{1}{1-X}$ then $\mu[f]$ is a finite state with the set of states
$Q_f=\{\mu[f]\circ\alpha[s]\,|\, s\in R\}$. 
\end{theorem}           

$\Box$  Let  $g=\sum\limits_{k=0}^\infty s_kX^k$  and since $f=\sum\limits_{k=0}^\infty X^k$ 
then

\begin{eqnarray*}
\mu[f](g) &=& fg=(1+\sigma(f)X)g=(1+fX)g \\
&=& (1+fX)(s_0+\sigma(g)X)\\
&=& s_0+\sigma(g)X+s_0fX+f\sigma(g)X^2\\
&=& s_0+(\sigma(g)+s_0f+f\sigma(g)X)X\\
&=& s_0+((1+fX)\sigma(g)+s_0f)X\\
&=& s_0+(f\sigma(y)+s_0f)X=s_0+f(s_0+\sigma(g))X\\
&=& s_0+\mu[f](s_0+\sigma(g))X\\
&=&s_0+\mu[f](\alpha[s_0](\sigma(g)))X
\end{eqnarray*}

We know
\[
\forall x\in R^\omega \quad \mu[f](rx)=\epsilon(a)+\mu_r[f](x),
\]
where $\mu_r[f]$ is the remainder of $\mu[f]$ for the letter $r\in R$.
Thus $\epsilon(s_0)=s_0$ and $\mu_{s_0}[f]=\mu[f]\circ \alpha[s_0]$.

Note if $s=0$ then $\mu[f]=\mu[f]\circ\alpha[s]$. Therefore we can apply \ref{lemma7.13}. Lemma for $Q=\{\mu[f]\circ\alpha[s]\,|\, s\in R\}$. This means we must prove 
\linebreak
$(\mu_s[f])_{s_9}\in Q$ for any pair $(s,s_0)\in R^2$.
\begin{eqnarray*}
\mu_s[f](g) &=& \mu[f](\alpha[s](g))=f(s+g)=(1+fX)(s+g)\\
&=& s+ sfX+(1+fX)g\\
&=& s+ sfX + s_0+f(s_0+\sigma(g))X\\
&=& s+s_0+f(s+s_0+\sigma(g))X\\
&=& s+s_0+\mu[f](s+s_0+\sigma(g))X\\
&=& s+s_0+\mu[f](\alpha[s+s_0](\sigma(g)))X
\end{eqnarray*}
Now we see that $(\mu_s[f])_{s_0}=\mu[f]\circ \alpha[s+s_0]\in Q$. Since $R$ is finite then $Q$ is finite.
\rule{2mm}{2mm}

\section{A Mealy machine}

This section concerns with $\displaystyle f(x)=\frac{1}{1-X}$ only.
The above theorem  tells us \cite{ski} that we can associate to $\mu[f]$ a Mealy
machine whose state functions are given by the states of $\mu[f]$. In other
words, define 

\[
\mathcal{M}_f = \langle Q_f, R,\circ, *\rangle
\] 
\begin{itemize}
	\item with states $Q_f = \{\alpha[s]\mu[f]\,|\, s\in R\}$ and
\item alphabet $R$,
\item $Q\times R \stackrel{\circ}{\longrightarrow}Q : \alpha[s]\mu[f]\circ r = \alpha[s+r]\mu[f]$,
\item
$Q\times A \stackrel{\ast }{\longrightarrow}A : \alpha[s]\mu[f]* r=s+r$. 
\end{itemize}

\begin{figure}[h]
\input{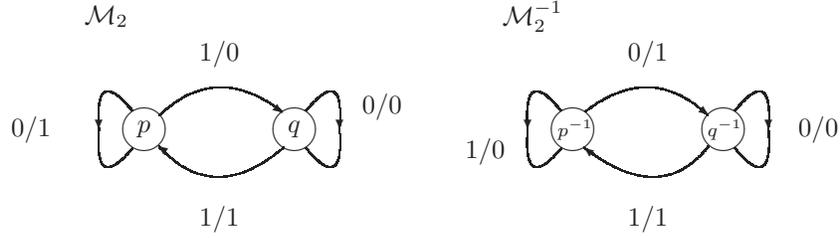}
\caption{A Mealy machine which generates the lamplighter group.}
\label{lamp2}
\end{figure}

Illustration when $R$ is the Galois field $GF(2)$ see in the figure \ref{lamp2} Here 
\[
\alpha[0]\mu[f]\mapsto q, \qquad \alpha[1]\mu[f] \mapsto p 
\]
and $\Gamma(\mathcal{M}_2)=\langle \bar q, \bar p\rangle=\langle \alpha[0]\mu[f],  \alpha[1]\mu[f]  \rangle$.

We know (see the proof of \ref{teor7.14}. Theorem)
\[
g\alpha[s]\mu[f]= s+r+\sigma(g)\alpha[s+r]\mu[f]X
\]
for any $g=r+\sum\limits_{k=1}^\infty s_kX^k\in R[[X]]$. This means that 
\[
g\alpha[s]\mu[f]= \alpha[s]\mu[f]* r+\sigma(g)\alpha[s+r]\mu[f]X
\]
 and the remainder $\big(\alpha[s]\mu[f]\big)_r=\alpha[s+r]\mu[f]$.
Now by \ref{prop5.5}. Proposition, \ref{prop5.8}. Proposition and \ref{teor3.13}. Theorem
we can conclude that
\[
g\alpha[s]\mu[f]=\alpha[s]\mu[f]*(r,s_1,s_2, \ldots, s_n, \ldots)
\]

\begin{lemma}\label{lem8.1}
If \quad
$K= \langle \alpha[sf], \mu[f] \,|\, s\in R\rangle$ \\
 is the group generated by $\{\alpha[sf], \mu[f]\,|\, s\in R\}$ and
\[
L= \langle \alpha[s]\mu[f] \,|\, s\in R\rangle
\]
 is the group generated by $\{\alpha[s]\mu[f]\,|\, s\in R\}$ then
$K=L$.
\end{lemma}

\begin{eqnarray*}
\Box\quad \mu[f] &=& \alpha[0]\mu[f]\in L, \\
\mu[f^{-1}] &=& \mu[f]^{-1}\in L: \qquad \ref{prop7.6}. {\textnormal{ Proposition,}}\\
\alpha[sf] &=& \alpha[fs]= \mu[f^{-1}]\alpha[s]\mu[f]\in L: \quad \ref{prop7.7}.  {\textnormal{ Proposition.}}
\end{eqnarray*}
Thus $K\subseteq L$.

\begin{eqnarray*}
\mu[f^{-1}] &=& \mu[f]^{-1}\in K\\
\alpha[s] &=& \alpha[sf^{-1}f]=\alpha[f^{-1}sf]=\mu[f]\alpha[sf]\mu[f^{-1}]\\
&=& \big(\mu[f]\alpha[sf]\big) \mu[f^{-1}]\in K,\\
\alpha[s]\mu[f] &\in& K.
\end{eqnarray*}
Thus $L \subseteq K$.

Hence $K \subseteq L \subseteq K$. This means $K=L$.
\rule{2mm}{2mm}

\begin{corollary}\label{cor8.2}
If $R=GF(2)$ then $L=\Gamma(\mathcal{M}_2)$.
\end{corollary}

\begin{lemma}
$\forall h_1 \in R[[X]] \; \forall h_2\in R[[X]] \quad \mu[h_1h_2]=\mu[h_1]\mu[h_2]$
\end{lemma}
$\Box$ Let $g\in R[[X]]$ then
\begin{eqnarray*}
g\mu[h_1h_2]&=& h_1h_2 g=h_2h_1g=\mu[h_2](h_1g)=\mu[h_2](\mu[h_1](g))\\
&=& g\mu[h_1]\mu[h_2]. \quad \rule{2mm}{2mm}
\end{eqnarray*}

\begin{lemma}\label{l8.4}
$\forall m\in \mathbb{Z}\; \mu[f^m]=\mu[f]^m$
\end{lemma}
\begin{eqnarray*}
\Box \quad \mu[f^2] &=&\mu[ff]=\mu[f]\mu[f]=\mu[f]^2,\\
\mu[f^{-2}] &=& \mu[f^{-1}f^{-1}]=\mu[f^{-1}]\mu[f^{-1}]=\mu[f]^{-1}\mu[f]^{-1}=\mu[f]^{-2}.
\end{eqnarray*}
The rest is inductiom.

(i) If $m>0$ then
\[
\mu[f^{m+1}]=\mu[f^mf]=\mu[f^m]\mu[f]=\mu[f]^m\mu[f]=\mu[f]^{m+1}.
\]

(ii) If $m<0$ then
\begin{eqnarray*}
\mu[f^{m-1}] &=& \mu[f^mf^{-1}]=\mu[f^m]\mu[f^{-1}]\\
&=& \mu[f]^m\mu[f]^{-1}=\mu[f]^{m-1}.\quad
\rule{2mm}{2mm}
\end{eqnarray*}

We shall restrict our attention further only on situation when $R=GF(2)$.
Let 
\[
\underset{\mathbb{Z}}{\bigoplus} \; GF(2) =\{\mathfrak{x}: \mathbb{Z} \to  GF(2)\,|\,\overset{\infty}{\forall}n\; \mathfrak{x}(n)=0 \}.
\]
This means that at most a finite number of $\mathfrak{x}(n)$ are not the 0 and we write
\[
\mathfrak{x}=(\mathfrak{x}_i)=\cdots \mathfrak{x}_{-n}\cdots \mathfrak{x}_{-1}\mathfrak{x}_0\mathfrak{x}_1\cdots\mathfrak{x}_n\cdots
\]
We denote by $GF^+(2)$ the aditive group of $GF(2)$.

\begin{definition} Let $(\mathfrak{x}_i), (\mathfrak{y_i})\in \underset{\mathbb{Z}}{\bigoplus} \; GF(2)$. Then we define $(\mathfrak{x_i})+(\mathfrak{y_i})=(\mathfrak{z_i})$
where, $\forall i\; \mathfrak{z}_i=\mathfrak{x}_i+\mathfrak{y}_i$, and  denote this group by $\underset{\mathbb{Z}}{\bigoplus} \; GF^+(2)$.
\end{definition} 

\begin{proposition}
The set $\{f^m\, |\, m \in \mathbb{Z}\}$ is linearly independent over $GF(2)$.
\end{proposition}

$\Box$ Since $f$ is invertible in $GF(2)[[X]]$ then $\forall m\in\mathbb{Z} \; f^m\ne 0$.
Let 
\[
f^{n_1}+f^{n_2}+\ldots +f^{n_k}=0
\]
where $n_1<n_2<\ldots < n_k$ then
\[
0=f^{|n_1|}(f^{n_1}+f^{n_2}+\ldots f^{n_k})=f^{\nu_1}+f^{\nu_2}+\ldots f^{\nu_k}.
\]
 Here $\forall i \; \nu_i=n_i+|n_1|$ and $0\le \nu_1<\nu_2<\ldots<\nu_k$.

\[
0=f^{-\nu_1}(f^{\nu_1}+f^{\nu_2}+\ldots f^{\nu_k})=1+f^{m_2}+f^{m_3}+\ldots + f^{m_k}.
\]
Here $\forall i \; m_i=\nu_i-\nu_1$.
\begin{align*}
(1-X)^{m_k} &+f^{m_2}(1-X)^{m_k}+\ldots + f^{m_k}(1-X)^{m_k}=0\\
(1-X)^{m_k} &+(1-X)^{m_k-m_2}+\ldots + (1-X)^{m_k-m_k}=0\\
(1-X)^{m_k} &+(1-X)^{m_k-m_2}+\ldots +(1-X)^{m_k-m_{k-1}}+ 1=0
\end{align*}
Let $p(X)=(1-X)^{m_k} +(1-X)^{m_k-m_2}+\ldots +(1-X)^{m_k-m_{k-1}}+ 1$ then $p(X)\in GF(2)[X]$, namely $p(X)$ is polinomial. If $X=1$ then 
\[
0=p(1)+1=1
\]
Contradiction! Since $GF(2)[X]$ is subring of $GF(2)[[X]]$ then we have gotten the contradiction in $GF(2)[[X]]$ too.
\rule{2mm}{2mm}

\begin{corollary}\label{cor8.6a}
$\langle\mu[f]\rangle\cong\mathbb{Z}$
\end{corollary}

$\Box$
Let $\mu[f]^n=\mu[f]^k$ then $\mu[f^{n-k}]=\mathbb{I}$, namely 
\[\mathbb{I}: GF(2)[[X]]\to GF(2)[[X]]:g\mapsto g\] is the identity map.
\[
1\mu[f^{n-k}]=f^{n-k}\ne 1= \mathbb{I}(1).
\]
Hence $\mu[f]$ has an infinite order.
\rule{2mm}{2mm}

\begin{corollary}\label{cor8.6}
If $(m_1,m_2,\ldots m_j)\ne(n_1,n_2,\ldots,n_i)$ then 
\[
\alpha[f^{m_1}]\alpha[f^{m_2}]\ldots\alpha[f^{m_j}]\ne\alpha[f^{n_1}]\alpha[f^{n_2}]\dots\alpha[f^{n_i}]
\]
\end{corollary}

\begin{eqnarray*}
\Box \quad 0\alpha[f^{m_1}]\alpha[f^{m_2}]\ldots\alpha[f^{m_j}] &=& f^{m_1}+f^{m_2}+\ldots+f^{m_j}\\
&\ne& f^{n_1}+f^{n_2}+\dots+f^{n_i}\\
&=&0\alpha[f^{n_1}]\alpha[f^{n_2}]\dots\alpha[f^{n_i}] \quad \rule{2mm}{2mm}
\end{eqnarray*}

\begin{proposition}\label{prop8.9}
Let $N=\langle \alpha[sf^m]\,|\, m\in \mathbb{Z}\rangle$ then
\[N \cong \underset{\mathbb{Z}}{\bigoplus} \; GF^+(2)\]
\end{proposition}

$\Box$ Since $g\alpha[f^n]\alpha[f^m]=g\alpha[f^m]\alpha[f^n]$ for any $g\in GF[[2]]$ then every element $\beta\in N$
can be represented for some $j$ as composition 
\[
\alpha[f^{m_1}]\alpha[f^{m_2}]\ldots\alpha[f^{m_j}]
\]  where
$m_1<m_2<\ldots<m_j$. By \ref{cor8.6} Corollary such representation is unique. Now we can conclude that the map
\begin{eqnarray*}
\varphi &:& N\to \underset{\mathbb{Z}}{\bigoplus} \; GF^+(2)\\
 &:& \alpha[f^{m_1}]\alpha[f^{m_2}]\ldots\alpha[f^{m_j}]\mapsto \cdots \mathfrak{x}_{-n}\cdots \mathfrak{x}_{-1}\mathfrak{x}_0\mathfrak{x}_1\cdots\mathfrak{x}_n\cdots
\end{eqnarray*}
where
\[
\mathfrak{x}_i=
\begin{cases}
1, & {\rm if}\; \exists \kappa \; i=m_\kappa;\\
0, & {\rm otherwise}
\end{cases}
\]
defines the isomorphism of groups $N,\underset{\mathbb{Z}}{\bigoplus} \; GF^+(2)$. \rule{2mm}{2mm}

\begin{lemma}
$N  = \langle \mu[f]^{-m}\alpha[sf]\mu[f]^m\,|\, m\in \mathbb{Z}\rangle$
\end{lemma}

$\Box$ Since $f$ is invertable then $f^m$ is invertable for every $m\in\mathbb{Z}$. We know (\ref{prop7.7}. Proposition)
that $\mu[f^{-m}]\alpha[h]\mu[f^m]=\alpha[f^mh]$ for any $h\in GF[[X]]$. Hence
\[ \mu[f]^{-m}\alpha[sf]\mu[f]^m=\mu[f^{-m}]\alpha[sf]\mu[f^m]=\alpha[f^msf]=\alpha[sf^{m+1}] \quad
\rule{2mm}{2mm}
\]

\begin{proposition}
$N\trianglelefteq \Gamma(\mathcal{M}_2)$. This means that $N$ is normal in $\Gamma(\mathcal{M}_2)$.
\end{proposition}

$\Box$ Let $\beta=\mu[f]^{-m}\alpha[s'f]\mu[f]^m\alpha[s]\mu[f]$ then $\beta\in N\alpha[s]\mu[f]$.
If $s'=0$ then $\beta=\alpha[s]\mu[f]\in\alpha[s]\mu[f]N$ otherwise $s'=1$.
\begin{eqnarray*}
\beta &=&\alpha[s]\alpha[-s]\mu[f]^{-m}\alpha[f]\mu[f]^m\alpha[s]\mu[f]\\
&=& \alpha[s]\mu[f]\mu[f]^{-1}\alpha[-s]\mu[f]^{-m}\alpha[f]\mu[f]^m\alpha[s]\mu[f]\\
&=& \alpha[s]\mu[f]\Bigl(\mu[f]^{-1}\alpha[-s]\mu[f]\Bigr)\\
&&
\Bigl(\mu[f]^{-m-1}\alpha[f]\mu[f]^{m+1}\Bigr)\Bigl(\mu[f]^{-1}\alpha[s]\mu[f]\Bigr)\\
&=& \alpha[s]\mu[f] \Bigl(\alpha[-sf]\Bigr)\Bigl(\mu[f]^{-m-1}\alpha[f]\mu[f]^{m+1}\Bigr) \alpha[sf]\Bigl)\\
&\in &\alpha[s]\mu[f] N
\end{eqnarray*}

If $\gamma\in N\alpha[s]\mu[f]$ then 
\[
\gamma=\Bigl(\mu[f]^{-m_1}\alpha[f]\mu[f]^{m_1}\Bigr)\ldots\Bigl(\mu[f]^{-m_i}\alpha[f]\mu[f]^{m_i}\Bigr)\alpha[s]\mu[f]
\]
for some $i$ where $\forall j\; m_j\in\mathbb{Z}$.
Hence
\begin{eqnarray*}
\gamma &=&\alpha[s]\mu[f]\mu[f]^{-1}\alpha[-s]\\
&&\Bigl(\mu[f]^{-m_1}\alpha[f]\mu[f]^{m_1}\Bigr)\ldots\Bigl(\mu[f]^{-m_i}\alpha[f]\mu[f]^{m_i}\Bigr)\alpha[s]\mu[f]\\
&=&\alpha[s]\mu[f] \Bigl(\mu[f]^{-1}\alpha[-s]\mu[f]\Bigr)\Bigl(\mu[f]^{-m_1-1}\alpha[f]\mu[f]^{m_1+1}\Bigr)\ldots\\
&\dots& \Bigl(\mu[f]^{-m_i-1}\alpha[f]\mu[f]^{m_i+1}\Bigr)\Bigl(\mu[f]^{-1}\alpha[s]\mu[f]\Bigr)\\
&=&\alpha[s]\mu[f] \Bigl(\alpha[-sf]\Bigr)\Bigl(\mu[f]^{-m_1-1}\alpha[f]\mu[f]^{m_1+1}\Bigr)\ldots\\
&\dots&\Bigl(\mu[f]^{-m_i-1}\alpha[f]\mu[f]^{m_i+1}\Bigr)\alpha[sf]\\
&\in &\alpha[s]\mu[f] N \quad \rule{2mm}{2mm}
\end{eqnarray*}

\section{A Semidirect product}

Let 
\begin{itemize}
	\item[$\bullet$] $\mathcal{N}$ and $\mathcal{H}$ be groups;
	\item[$\bullet$] $\varphi: \mathcal{H} \to \mathfrak{A}ut(\mathcal{N})$ ---  a group homomorphism, where 
	$\mathfrak{A}ut(\mathcal{N})$ --- the automorphism group of the group $\mathcal{N}$.
\end{itemize}
We shall write $\varphi^h\in\mathfrak{A}ut(\mathcal{N})$ insted of $\varphi(h)$.
Then define $\mathcal{G}=\mathcal{N}\rtimes_\varphi \mathcal{H}=\mathcal{N}\rtimes \mathcal{H}$ to be the set
$\mathcal{N}\times \mathcal{H}$ with the multiplication defined by
\[
(n_1,h_1)(n_2,h_2)=(n_1\varphi^{h_1}(n_2),h_1h_2).
\]

\begin{definition}\label{def9.1}
Let $\mathcal{G}$ be a group with subgroups $\mathcal{N}$ and $\mathcal{H}$ such that 
\[
\mathcal{N}\cap \mathcal{H}=\{e\} \quad {\rm and} \quad \mathcal{N}\mathcal{H}=\mathcal{G}.
\]
If $\mathcal{N}$ is normal (but not necessarily $\mathcal{H}$), then we say that $\mathcal{G}$ is the semidirect product of $\mathcal{N}$ and $\mathcal{H}$. Here $e$ is the neutral element of $\mathcal{G}$.
\end{definition}

We recall some standard results (see for example in \cite{Con}).

\begin{lemma}
Let $h\in\mathcal{G}$. If $\mathcal{N}\trianglelefteq \mathcal{G}$ then $\varphi^h\in \mathfrak{A}ut(\mathcal{N})$ where \\ $\varphi^h: x\mapsto hx h^{-1}$.
\end{lemma}

$\Box$ (i) $\varphi^h(xy)=hxyh^{-1}=hxh^{-1}hyh^{-1}=\varphi^h(x)\varphi^h(y)$. Thus $\varphi^h$ is a homomorphism.

(ii) If $\varphi^h(x)=\varphi^h(y)$ then $hxh^{-1}=hyh^{-1}$. Hence $x=y$, therefore $\varphi^h$ is an injection.

(iii) Since $\mathcal{N}\trianglelefteq \mathcal{G}$ then $h^{-1}yh\in\mathcal{N}$ for every $y\in\mathcal{N}$.
Hence $\varphi^h(h^{-1}yh)=hh^{-1}yhh^{-1}=y$. Thus $\varphi^h$ is a surjection.
\rule{2mm}{2mm}

\begin{theorem}\label{teor9.3}
Let $\mathcal{G}$ be the semidirect product of $\mathcal{N}$ and $\mathcal{H}$. If 
\[
\varphi: \mathcal{H} \to \mathfrak{A}ut(\mathcal{N}): h\mapsto \varphi^h
\]
is defined by
\[
\varphi^h:\mathcal{N} \to \mathcal{N}: n\mapsto hnh^{-1}
\]
then  $\varphi$ is a homomorphism and 
\[
f: \mathcal{N}\rtimes_\varphi \mathcal{H} \to \mathcal{G}: (n,h)\mapsto nh
\]
is an isomorphism.
\end{theorem}

 \begin{eqnarray*}
\Box \; {\rm (i)}\quad \varphi^{kh}(n) &=& khn(kh)^{-1}=khnh^{-1}k^{-1},\\
\varphi^k\circ\varphi^h(n) &=& \varphi^k(\varphi^h(n))=\varphi^k(hnh^{-1})=khnh^{-1}k^{-1}.
\end{eqnarray*}
Thus $\varphi^{kh}=\varphi^k\circ \varphi^h$ or $\varphi(kh)=\varphi(k)\circ\varphi(h)$. Hence
\begin{center}
$\varphi: \mathcal{H} \to \mathfrak{A}ut(\mathcal{N}): h\mapsto \varphi(h)$
\end{center} is a homomorphism.

(ii) Since $\mathcal{G}$ is the semidirect product of $\mathcal{N}$ and $\mathcal{H}$ then $\mathcal{G}=\mathcal{N}\mathcal{H}$. Therefore the map
\[
f: \mathcal{N}\rtimes_\varphi \mathcal{H} \to \mathcal{G}: (n,h)\mapsto nh
\] is surjective.

(iii) Let $f(n_1,h_1)=f(n_2,h_2)$. This means that $n_1h_1=n_2h_2$ or $n_2^{-1}n_1=h_2h_1^{-1}$.
Since $\mathcal{G}$ is the semidirect product of $\mathcal{N}$ and $\mathcal{H}$ then $\mathcal{N}\cap\mathcal{H}=\{e\}$.
Hence $n_2^{-1}n_1=e=h_2h_1^{-1}$. Thus $n_1=n_2$ and $h_1=h_2$. Therefore $f$ is injective.

\begin{eqnarray*}
{\rm (iv)} \quad f((n_1,h_1)(n_2,h_2)) &=& f(n_1\varphi^{h_1}(n_2),h_1h_2)=n_1\varphi^{h_1}(n_2)h_1h_2\\
&=& n_1(h_1n_2h_1^{-1})h_1h_2=n_1h_1n_2h_2\\
&=& f(n_1,h_1)f(n_2,h_2)
\end{eqnarray*}
Hence $f$ is a homomorphism.

Now we can conclude that $f$ is an isomorphism.
\rule{2mm}{2mm}

\begin{lemma}
If $\phi\in \mathfrak{A}ut(\mathcal{N}_1)$ and $f:\mathcal{N}_1 \to \mathcal{N}_2$ is an isomorphism then
$f\circ \phi\circ f^{-1}\in \mathfrak{A}ut(\mathcal{N}_2)$.
\end{lemma}

$\Box$
Since $f,\phi$ are bijective then $f\circ \phi\circ f^{-1}$ is bijective and $f\circ \phi \circ f^{-1}$ is a homomorphism because
this is composition of homomorphisms.
\rule{2mm}{2mm}

\begin{proposition}\label{prop9.5}
Let $\mathcal{N}_1, \mathcal{N}_2, \mathcal{H}_1, \mathcal{H}_2$ be groups. If 
\[\mathcal{N}_1\cong\mathcal{N}_2 \quad  {\rm and} \quad \mathcal{H}_1\cong\mathcal{H}_2\]
 then for every homomorphism $\varphi:\mathcal{H}_1 \to \mathfrak{A}ut(\mathcal{N}_1)$ there exists the homomorphism 
 $\psi:\mathcal{H}_2\to \mathfrak{A}ut(\mathcal{N}_2)$ such that 
\[
\mathcal{N}_1\rtimes_\varphi\mathcal{H}_1\cong\mathcal{N}_2\rtimes_\psi\mathcal{H}_2.
\]
\end{proposition}

$\Box$ Let $\psi^{h'}=f\circ \varphi^{g^{-1}(h')}\circ f^{-1}$ and $F:(n,h)\mapsto (f(n),g(h))$ where
$f:\mathcal{N}_1\to \mathcal{N}_2$, $g:\mathcal{H}_1\to \mathcal{H}_2$ are isomorphisms then
\begin{eqnarray*}
F((n_1,h_1)(n_2,h_2)) &=& F(n_1\varphi^{h_1}(n_2),h_1h_2)=(f(n_1\varphi^{h_1}(n_2)),g(h_1h_2))\\
&=&(f(n_1)f(\varphi^{h_1}(n_2)),g(h_1h_2))\\
&=&(f(n_1)(f\circ\varphi^{h_1}(n_2)),g(h_1h_2))\\
F(n_1,h_1)F(n_2,h_2) &=& (f(n_1),g(h_1))(f(n_2),g(h_2))\\
&=& (f(n_1)\psi^{g(h_1)}(f(n_2)),g(h_1)g(h_2))\\
&=& (f(n_1)(f\circ\varphi^{g^{-1}\circ g(h_1)}\circ f^{-1}\circ f(n_2)),g(h_1h_2))\\
&=& (f(n_1)(f\circ \varphi^{h_1}(n_2)),g(h_1h_2))
\end{eqnarray*}
This means that $F:\mathcal{N}_1\rtimes_\varphi\mathcal{H}_1\to\mathcal{N}_2\rtimes_\psi\mathcal{H}_2$ is a homomorphism.
Since $F$ is bijective than $F$ is an isomorphism.
\rule{2mm}{2mm}
\medskip

Now we return to $\Gamma(\mathcal{M}_2)$.

\begin{proposition}
$\langle \mu[f]\rangle$ is not a normal subgroup of $\Gamma(\mathcal{M}_2)$.
\end{proposition}

$\Box$ Let $a=\alpha[f]\mu[f]\alpha[f]$ and $g\in GF(2)[[X]]$ then $a\mu[f]\in a\langle \mu[f] \rangle$ and
\begin{eqnarray*}
 ga\mu[f] &=& g\alpha[f]\mu[f]\alpha[f]\mu[f]=((g+f)f+f)f
= gf^2+f^3+f^2,\\
g\mu[f^k]a &=& g\mu[f]\alpha[f]\mu[f]\alpha[f]=(gf^k+f)f+f=gf^{k+1}+f^2+f\\
&\ne&gf^2+f^3+f^2.
\end{eqnarray*}
If $k=1$ then $g\mu[f^k]a=gf^2+f^2+f$ but even in this situation\\ $gf^2+f^2+f^3\ne gf^2+f^2+f$.
Hence $a\langle\mu[f]\rangle\ne \langle\mu[f]\rangle a$.
\rule{2mm}{2mm}

\begin{theorem} \label{t9.7}
$\Gamma(\mathcal{M}_2)\cong \underset{\mathbb{Z}}{\bigoplus} \; GF^+(2)\rtimes \mathbb{Z}$
\end{theorem}

$\Box$ (i) If $x\in N\cap\langle \mu[f]\rangle$ then $x=\mathbb{I}$ or $x=\alpha[f^m]$ for some $m$ and $x=\mu[f^n]$
for some $n$. Hence
\[
1+f^m=1\alpha[f^m]=1x=1\mu[f^k]=f^k
\]
Contradiction! Thus $N\cap\langle \mu[f]\rangle=\{\mathbb{I}\}$.

(ii) Let $x\in\Gamma(\mathcal{M}_2)$ then $x=a_1a_2\ldots a_n$ for some $n$ where $a_i\in \{\alpha[f],\mu[f]\}$
(see \ref{lem8.1}. Lemma and \ref{cor8.2}. Corollary). Thus by \ref{prop8.9}. Proposition $\forall i \; a_i\in N\cup \langle \mu[f]\rangle$. If $a_i\in\langle \mu[f]\rangle $ but $a_{i+1}\in N$ then there exists $\alpha_i\in N$ such that
\[
a_ia_{i+1}=\alpha_ia_i
\]
because $N\trianglelefteq \Gamma(\mathcal{M}_2)$. Hence $x\in N\langle \mu[f]\rangle$ and $\Gamma(\mathcal{M}_2)=N\langle \mu[f]\rangle$

(iii) Since $N\cap\langle \mu[f]\rangle=\{\mathbb{I}\}$, $\Gamma(\mathcal{M}_2)=N\langle \mu[f]\rangle$ and 
$N\trianglelefteq \Gamma(\mathcal{M}_2)$ then by \ref{def9.1}. Definition $\Gamma(\mathcal{M}_2)$ is the semidirect product of $N$ and $\langle \mu[f]\rangle$. Hence
\begin{eqnarray*}
\Gamma(\mathcal{M}_1) &\cong& N\rtimes \langle \mu \rangle \qquad {\rm \ref{teor9.3}. \; Theorem}\\
N\rtimes \langle \mu \rangle &\cong& \underset{\mathbb{Z}}{\bigoplus} \; GF^+(2)\rtimes \mathbb{Z} 
\end{eqnarray*}
by  \ref{cor8.6a}. Corollary, \ref{prop8.9}. Proposition and  \ref{prop9.5}. Proposition.
Therefore $\Gamma(\mathcal{M}_2)\cong \underset{\mathbb{Z}}{\bigoplus} \; GF^+(2)\rtimes \mathbb{Z}$. 
\rule{2mm}{2mm}

\section{The Lamplighter group}

We recall basic results (for more information see \cite{Eck}).

\begin{figure}[h]
\includegraphics[trim=2cm 19.3cm 0cm 4.3cm, clip ]{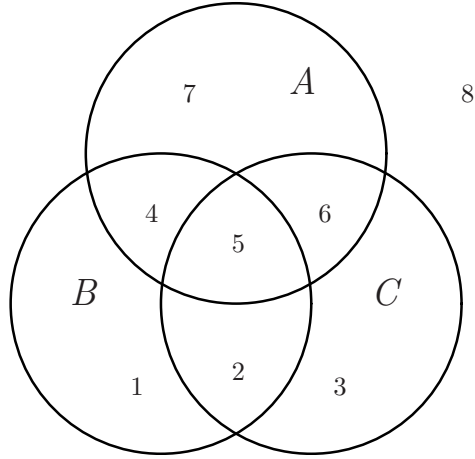}\\
\caption{$A=\{4,5,6,7\}$,
$B=\{1,2,4,5\}$, $C=\{2,3,5,6\}$.}
\label{z89e}
\end{figure}

\begin{proposition}\label{a5.18.1}
Let $X$ be a set and 
$\mathfrak{P}(X)=\{A\,|\, A\subseteq X\}$. An algebra $\langle \mathfrak{P}(X), A\bigtriangleup B\rangle$ is a commutative group.
\end{proposition}

$\Box$ (i) Associativity. See the  Venn diagramm (\ref{z89e} fig.):
\begin{eqnarray*}
A&=&\{4,5,5,7\},\\
B&=&\{1,2,4,5 \},\\
C&=&\{ 2,3,5,6 \}
\end{eqnarray*}
\begin{alignat*}{3}
A\setminus B &=\{ 6,7\}, &\quad  B\setminus A &=\{1,2 \}, &\quad A \bigtriangleup B &=\{1,2, 6, 7 \},\\
B\setminus C &= \{1,4 \}, & C \setminus B &= \{3,6 \}, & B\bigtriangleup C &=\{1, 3, 4, 6 \}
\end{alignat*}
\begin{alignat*}{1}
 (A \bigtriangleup B )\setminus C &=\{1,2,6,7 \}\setminus \{2,3,5,6 \} =\{1,7 \},\\
C\setminus (A \bigtriangleup B) &=\{ 2,3,5,7\}\setminus \{ 1,2,6,7\}= \{3,5 \},\\
(A \bigtriangleup B)\bigtriangleup C &=\{1,3,5,7 \},\\
A \setminus (B\bigtriangleup C) &= \{4,5,6,7 \}\setminus \{1,3,4,6 \}=\{5, 7 \},\\
(B\bigtriangleup C)\setminus A &= \{1,3,4,6 \}\setminus \{4,5,6,7 \}=\{1,3 \},\\
A \bigtriangleup (B\bigtriangleup C) &= \{1,3,5,7 \}
\end{alignat*}
Thus $A\bigtriangleup (B\bigtriangleup C)=(A\bigtriangleup B)\bigtriangleup C$.

(ii) $\emptyset$ is a neutral element of the semigroup. 

(iii) $A$ is the inverse element of a set $A$.

(iv) $A\bigtriangleup B=(A\setminus B)\cup (B\setminus A)= (B\setminus A)\cup (A\setminus B)=B\bigtriangleup A$.
\rule{2mm}{2mm}

\begin{corollary}
$A\bigtriangleup B=A\bigtriangleup C \Rightarrow B=C$.
\end{corollary}

Let $L_2=\{(S,x)\,|\,S\subset\mathbb{Z}\wedge|S|<\aleph_0\wedge x\in\mathbb{Z}\}$.

Let  $l_1=(S,x), l_2=(T,y)$ are elements of $L_2$ then we define 
\[l_1l_2=(S_y\bigtriangleup T, x+y),\]
 where $S_y=\{s+y\,|\,s\in S\}$. Hence we have a groupoid  $L_2$.

\begin{lemma}
$(S\wedge T)_y=S_y\wedge T_y$ and $(S\setminus T)_y=S_y\setminus T_y$
\end{lemma}

$\Box$ (i) Let $z\in(S\wedge T)_y$ then $z-y\in S\wedge T$ thence $z-y\in S$ and $z-y\in T$. Hence: $z\in S_y$ and $z\in T_y$ therefore $z\in S_y\wedge T_y$. Thus
 $(S\wedge T)_y\subseteq S_y\wedge T_y$.

(ii) Let $z\in S_y\wedge T_y$ then $z\in S_y$ and $z\in T_y$ thence $z-y\in S$ and $z-y\in T$. Hence: $z-y\in S\wedge T$ therefore $z\in (S\wedge T)_y$.
Thus $S_y\wedge T_y\subseteq  (S\wedge T)_y$.

 $(S\wedge T)_y=S_y\wedge T_y$ from (i) and (ii).

(iii) Let  $z\in(S\setminus T)_y$ then $z-y\in S\setminus T$ thence $z-y\in S$ and $z-y\not\in T$. hence: $z\in S_y$ and $z\not\in T_y$ therefore $z\in S_y\setminus T_y$.
Thus $(S\setminus T)_y\subseteq S_y\setminus T_y$.

(iv)  Let $z\in S_y\setminus T_y$ then $z\in S_y$ and $z\not\in T_y$ thence $z-y\in S$ and $z-y\not\in T$. Hence: $z-y\in S\setminus T$ therefore $z\in (S\setminus T)_y$.
Thus $S_y\setminus T_y\subseteq  (S\setminus T)_y$.

$(S\setminus T)_y=S_y\setminus T_y$ from (iii) and (iv).
\rule{2mm}{2mm}

\begin{corollary}\label{s5.18.4}
$(S\bigtriangleup T)_y=S_y\bigtriangleup T_y$
\end{corollary}

\begin{eqnarray*}\Box \quad
(S\bigtriangleup T)_y&=&((S\setminus T)\wedge (T\setminus S))_y=(S\setminus T)_y\wedge(T\setminus S)_y\\
&=& (S_y\setminus T_y)\wedge(T_y\setminus S_y)\\
&=& S_y \bigtriangleup T_y \quad \rule{2mm}{2mm}
\end{eqnarray*}

\begin{lemma}\label{l5.18.5}
$(S_x)_y=S_{x+y}$
\end{lemma}

$\Box$ $z\in(S_x)_y$ if and only if then $z-y\in S_x$. Hence $z-y\in S_x$ if and only if then $(z-y)-x\in S$. 
Successively $z\in S_{x+y}$ if and only if then $z-(x+y)\in S$.
Since $(z-y)-x=z-(x+y)$ then we have equivalence
\[
z\in(S_x)_y\Leftrightarrow z\in S_{x+y}.
\]
Thus $(S_x)_y=S_{x+y}$ \rule{2mm}{2mm}

\begin{proposition}
$L_2$ is a group.
\end{proposition}

$\Box$ (i) Let $l_1=(S,x), l_2=(T,y), l_3=(R,z)$ are elements of $L_2$ then 
\begin{alignat*}{2}
l_1l_2 &=(S_y\bigtriangleup T, x+y), & \quad \textnormal{where} \quad S_y&= \{s+y\,|\,s\in S\},\\
l_2l_3 &=(T_z\bigtriangleup R, y+z), & \quad \textnormal{where} \quad T_z&= \{t+z\,|\,t\in T\}.
\end{alignat*}
\[
(l_1l_2)l_3 =(S_y\bigtriangleup T,x+y)(R,z)=((S_y\bigtriangleup T)_z\bigtriangleup R,x+y+z),
\]
where $(S_y\bigtriangleup T)_z=\{q+z\,|\,q\in S_y\bigtriangleup T\}$.
\[
l_1(l_2l_3)=(S,x)(T_z\bigtriangleup R, y+z)=(S_{y+z}\bigtriangleup T_z\bigtriangleup R, x+y+z),
\]
where $S_{y+z}=\{\sigma+y+z\,|\,\sigma\in S\}$.

(ii)
\[ 
(S_y\bigtriangleup T)_z\underset{C\ref{s5.18.4}}{=}(S_y)_z\bigtriangleup T_z\underset{L\ref{l5.18.5}}{=}S_{y+z}\bigtriangleup T_z
\]

(iii) 
\[
l_1(l_2l_3)\underset{\rm (i)}{=}(S_{y+z}\bigtriangleup T_z\bigtriangleup R,x+y+z)\underset{\rm (ii)}{=}((S_y\bigtriangleup T)_z\bigtriangleup R,x+y+z)\underset{\rm (i)}{=}(l_1l_2)l_3  
\]
Hence $L_2$ is a semigroup.

(iv) Let $n\in\mathbb{Z}$ then $S_n=\{s+n\,|\,s\in S\}$. 
$(\emptyset,0)\in L_2$ and $(\emptyset,0)(S,x)=(\emptyset _x\bigtriangleup S, 0+x)=(\emptyset\bigtriangleup S,x)=(S,x)$. Similarly
$(S,x)(\emptyset,0)=(S_0\bigtriangleup \emptyset,x+0)=(S\bigtriangleup \emptyset,x)=(S,x)$. Thus $L_2$ is a monoid.

(v) $(S_{-x},-x)\in L_2$ where $S_{-x}=\{s-x\,|\,s\in S\}$ and 
\[
(S,x)(S_{-x},-x)=(S_{-x}\bigtriangleup S_{-x},0)=(\emptyset,0).
\] 
Thus $(S_{-x},-x)$ is an  inverse element of $(S,x)$. 
\rule{2mm}{2mm}

\begin{proposition}
$L_2$ is not a commutative group.
\end{proposition}

$\Box$ Let $l_0=(\{0\},0), l_1=(\{1\},1)$ then $\{0\}_1=\{1\}, \{1\}_0=\{1\}$.  
\begin{eqnarray*}
l_0l_1 &=&(\{0\}_1\bigtriangleup \{1\}, 0+1)=(\{1\}\bigtriangleup \{1\},1)=(\emptyset,1),\\
l_1l_0 &=&  (\{1\}_0\bigtriangleup \{0\}, 1+0)=(\{1\}\bigtriangleup \{0\},1)=(\{0,1\},1).
\end{eqnarray*}
Hence $l_0l_1\ne l_1l_0$.
\rule{2mm}{2mm}

\begin{proposition}
$L_{\bar 2}$ is group where 

\centering{ $(S,x)\bar\circ(T,y)=(S_{-y}\bigtriangleup T,x+y)$.}
\end{proposition}

$\Box$ (i) Let $l_1=(S,x), l_2=(T,y), l_3=(R,z)$ are elements of $L_{\bar 2}$ then 
\begin{alignat*}{2}
l_1\bar \circ l_2 &=(S_{-y}\bigtriangleup T, x+y), & \quad \textnormal{where} \quad S_{-y}&= \{s-y\,|\,s\in S\},\\
l_2\bar \circ l_3 &=(T_{-z}\bigtriangleup R, y+z), & \quad \textnormal{where} \quad T_{-z}&= \{t-z\,|\,t\in T\}.
\end{alignat*}
\[
(l_1\bar \circ l_2)\bar \circ l_3 =(S_{-y}\bigtriangleup T,x+y)\bar \circ (R,z)=((S_{-y}\bigtriangleup T)_{-z}\bigtriangleup R,x+y+z),
\]
where $(S_{-y}\bigtriangleup T)_{-z}=\{q-z\,|\,q\in S_{-y}\bigtriangleup T\}$.
\[
l_1\bar \circ (l_2\bar \circ l_3)=(S,x)\bar \circ (T_{-z}\bigtriangleup R, y+z)=(S_{-y-z}\bigtriangleup T_{-z}\bigtriangleup R, x+y+z),
\]
where $S_{-y-z}=\{\sigma-y-z\,|\,\sigma\in S\}$.

(ii)
\[ 
(S_{-y}\bigtriangleup T)_{-z}\underset{C\ref{s5.18.4}}{=}(S_{-y})_{-z}\bigtriangleup T_{-z}\underset{L\ref{l5.18.5}}{=}S_{-y-z}\bigtriangleup T_{-z}
\]

(iii) 
\begin{eqnarray*}
l_1\bar \circ (l_2\bar \circ l_3)&\underset{\rm (i)}{=}&(S_{-y-z}\bigtriangleup T_{-z}\bigtriangleup R,x+y+z)\\
&\underset{\rm (ii)}{=}&((S_{-y}\bigtriangleup T)_{-z}\bigtriangleup R,x+y+z)\underset{\rm (i)}{=}(l_1\bar \circ l_2)\bar \circ l_3 
\end{eqnarray*}
Hence $L_{\bar 2}$ is a semigroup.

(iv) Let $x\in\mathbb{Z}$ then $S_x=\{s+x\,|\,s\in S\}$. 
$(\emptyset,0)\in L_{\bar 2}$ and $(\emptyset,0)\bar \circ(S,x)=(\emptyset _{-x}\bigtriangleup S, 0+x)=(\emptyset\bigtriangleup S,x)=(S,x)$. Similarly
$(S,x)\bar \circ (\emptyset,0)=(S_{-0}\bigtriangleup \emptyset,x+0)=(S\bigtriangleup \emptyset,x)=(S,x)$. Thus $L_{\bar 2}$ is a monoid.

(v) 
\[
(S,x)\bar \circ (S_{x},-x)=(S_{x}\bigtriangleup S_{x},0)=(\emptyset,0).
\] 
Thus $(S_{x},-x)$ is an  inverse element of $(S,x)$. 
\rule{2mm}{2mm}

\begin{lemma}
$(-S)_{-y}=-(S_y)$ where $-S=\{-s\,|\, s\in S\}$.
\end{lemma}
\begin{eqnarray*} \Box \;
(-S)_{-y} &=& \{-s\,|\, s\in S\}_{-y}=\{-s-y\,|\, s\in S\},\\
-(S_y) &=& -(\{s\,|\, s\in S\}_y)=-\{s+y\,|\, s\in S\}=\{-s-y\,|\, s\in S\} \; \rule{2mm}{2mm}
\end{eqnarray*}

\begin{proposition}
$L_2\cong L_{\bar 2}$
\end{proposition}

$\Box$ Let $f:L_2\to L_{\bar 2}: (S,x)\mapsto (-S,y)$ then $f$ is a bijection.
\begin{eqnarray*}
f((S,x)(T,y)) &=& f(S_y\bigtriangleup T, x+y)=(-(S_y\bigtriangleup T),x+y),\\
f(S,x)\bar \circ f(T,y) &=& (-S,x)\bar \circ (-T,y)=((-S)_{-y}\bigtriangleup (-T),x+y)\\
&=&(-(S_y)\bigtriangleup (-T),x+y)\\
&=& (-(S_y\bigtriangleup T),x+y) \quad \rule{2mm}{2mm}
\end{eqnarray*}

\medskip
Let
\[
L_2'= \{ ((\mathfrak{x}_i),n) \,|\, (\mathfrak{x}_i)\in \underset{\mathbb{Z}}{\bigoplus} \; GF^+(2) \;\wedge\; n\in\mathbb{Z}\}
\]
 Let $l_1=((\mathfrak{x}_i),n), l_2=((\mathfrak{y}_i),m)$ then we define $l_1l_2=((\mathfrak{z}_i),n+m)$ where $\mathfrak{z}_i=\mathfrak{x}_{i-m}\oplus \mathfrak{y}_i$ where $\oplus$ is addition in $GF^+(2)$.

\begin{proposition}
$L'_2$ is a group.
\end{proposition}

$\Box$ (i) Let $l_1=((\mathfrak{x}_i),n), l_2=((\mathfrak{y}_i),m)$ then $l_1l_2=((\mathfrak{z}_i),n+m)$ where $\mathfrak{z}_i= \mathfrak{x}_{i-m}\oplus \mathfrak{y}_i$.
If $\nu=|\{\mathfrak{x}_i\,|\, \mathfrak{x}_i=1\}|$ and  $\mu=|\{\mathfrak{y}_i\,|\, \mathfrak{y}_i=1\}|$ then
\[
|\{\mathfrak{z}_i\,|\, \mathfrak{z}_i=1\}|=|\{\mathfrak{x}_{i-m}\oplus \mathfrak{y}_i\,|\, \mathfrak{x}_{i-m}\oplus \mathfrak{y}_i=1\}|\le \nu+\mu.
\]
Thus $l_1l_2\in L'_2$, id est, $L'_2$ is a groupoid.

(ii) Let $l_3=((\mathfrak{k}_i), k)$ then 
\[
(l_1l_2)l_3 = ((\mathfrak{z}_i), n+m)((\mathfrak{k}_i),k)=((\mathfrak{t}_i),n+m+k),
\]
where $\mathfrak{z}_i=\mathfrak{x}_{i-m}\oplus \mathfrak{y}_i$ un $\mathfrak{t}_i=\mathfrak{z}_{i-k}\oplus \mathfrak{k}_i=\mathfrak{x}_{i-m-k}\oplus \mathfrak{y}_{i-k}\oplus \mathfrak{k}_i$.

On the other hand
\[
l_1(l_2l_3)=((\mathfrak{x}_i),n)((\mathfrak{v}_i),m+k)=((\mathfrak{w}_i),n+m+k),
\]
where $\mathfrak{v}_i=\mathfrak{y}_{i-k}\oplus \mathfrak{k}_i$ and $\mathfrak{w}_i=\mathfrak{x}_{i-m-k}\oplus\mathfrak{v}_i=\mathfrak{x}_{i-m-k}\oplus \mathfrak{y}_{i-k}\oplus \mathfrak{k}_i$. 
Therefore $(\mathfrak{t}_i)=(\mathfrak{w}_i)$. Hence $(l_1l_2)l_3=l_1(l_2l_3)$.
Thus  $L'_2$ is a semigroup.

(iii) Let $(0)=(\mathfrak{o}_i)\in \underset{\mathbb{Z}}{\bigoplus} \; GF^+(2)$ where $\forall i\; \mathfrak{0}_i=0$ then
\[
((0),0)((\mathfrak{x}_i),n)=((\mathfrak{x}_i),n)=((\mathfrak{x}_i),n)((0),0). 
\]
Thus $L'_2$ is a monoid.

(v) Let $l_1^{-1}=((\mathfrak{t}_i),-n)$ where $\mathfrak{t}_i= \mathfrak{x}_{i+n}$ then 
\[
l_1l_1^{-1}=((\mathfrak{x}_i),n)((\mathfrak{t}_i),-n)=((\mathfrak{o}_i),0),
\]
where $\mathfrak{o}_i=\mathfrak{x}_{i+n}\oplus \mathfrak{t}_i=\mathfrak{x}_{i+n}\oplus \mathfrak{x}_{i+n}=0$. Therefore $((\mathfrak{o}_i),0)=((0),0)$, id est, $((\mathfrak{o}_i),0)$ is a neutral element. 
Thus $l_1^{-1}$ is an  inverse element of $l_1$. 
\rule{2mm}{2mm}

\begin{proposition}
$L'_2\cong L_2$
\end{proposition}

$\Box$ (i) At first we define a map $\varphi: L'_2\to L_2$. Let $l_1=((\mathfrak{x}_i),n)\in L'_2$ then $\varphi((\mathfrak{x}_i),n)= (S,n)$ where $S=\{s\,|\,x_s=1\}$.

(ii)  Let $l_2=((\mathfrak{y}_i),m)$ and $\varphi(l_2)=(T,m)$ where $T= \{t\,|\,\mathfrak{y}_t=1\}$ then
$l_1l_2=((\mathfrak{z}_i),n+m)$ where $\mathfrak{z}_i=\mathfrak{x}_{i-m} \oplus \mathfrak{y}_i$ and $\varphi(l_1l_2)=(R,n+m)$ where
\[
R= \{\varrho\,|\, \mathfrak{z}_\varrho=1\}=\{\varrho\,|\, \mathfrak{x}_{\varrho-m}\oplus \mathfrak{y}_\varrho=1\}.
\]
On the other hand $\varphi(l_1)\varphi(l_2)=(S,n)(T,m)=(S_m\bigtriangleup T, n+m)$ where $S_m=\{i+m\,|\, i\in S\}$.
Thus $\varrho\in S_m\Leftrightarrow  \varrho-m\in S$.

We should like to show that   $R=S_m\bigtriangleup T$. Since $\varrho\in R\Leftrightarrow \mathfrak{x}_{\varrho-m} \oplus \mathfrak{y}_\varrho=1$ then 
\[
\mathfrak{x}_{\varrho-m}=0\; \wedge\;\mathfrak{y}_\varrho=1 \quad  {\rm or} \quad \mathfrak{x}_{\varrho-m}=1\; \wedge\;\mathfrak{y}_\varrho=0
\]
We know 
\begin{eqnarray*}
\mathfrak{x}_{\varrho-m}=1 &\Leftrightarrow & \varrho-m\in S \Leftrightarrow \varrho\in S_m,\\
\mathfrak{y}_\varrho =1 &\Leftrightarrow & \varrho \in T.
\end{eqnarray*}

The analysis of the situation is summarized in the table:

\begin{center}
\begin{tabular}{|c|c|l|l|} \cline{1-2}
$\mathfrak{x}_{\varrho-m}$ & $\mathfrak{y}_\varrho$ & \multicolumn{2}{l}{} \\ \hline  \hline
0 & 1 & $\varrho-m\not\in S$ & $\varrho\in T$ \\ 
 & & $\varrho\not\in S_m$ & \\ \cline{3-4}
 & & \multicolumn{2}{l|}{$\varrho\in T\setminus S_m\subseteq S_m\bigtriangleup T$} \\ \hline
1 & 0 & $\varrho-m\in S$ & $\varrho\not\in T$ \\ 
 & & $\varrho\in S_m$ & \\ \cline{3-4}
 & & \multicolumn{2}{l|}{$\varrho\in S_m\setminus T\subseteq S_m\bigtriangleup T$} \\ \hline
\end{tabular} \end{center}

Hence $\varrho\in R\Leftrightarrow \varrho\in S_m\bigtriangleup T$. Thence
\[
\varphi(l_1l_2)=(R,n+m)=(S_m\bigtriangleup T,n+m)=\varphi(l_1)\varphi(l_2),  
\]
id est,
$\varphi$ is a homomorphism of the groups $L'_2, L_2$.

(iii) Let $l_1=((\mathfrak{x}_i),n)$,  $l_2=((\mathfrak{y}_i),m)$ are elements of the set  $L'_2$ and  $l_1\ne l_2$. Let $\varphi(l_1)= (S,n)$ where $S=\{s\,|\,x_s=1\}$ and 
$\varphi(l_2)=(T,m)$ where $T= \{t\,|\,y_t=1\}$. If $n\ne m$ then $\varphi(l_1)\ne\varphi(l_2)$.

If $n=m$ and $\varphi(l_1)\ne \varphi(l_2)$ then $\exists k \; x_k\ne y_k$. Since $\mathfrak{x}_k,\mathfrak{y}_k$ are elements of $\{0,1\}$ then 
\begin{eqnarray*}
\mathfrak{x}_k=0 &\wedge& \mathfrak{y}_k=1,\\
{\rm or}\\
\mathfrak{x}_k=1 &\wedge& \mathfrak{y}_k=0.
\end{eqnarray*}
\begin{itemize}
\item If $\mathfrak{x}_k=0 \wedge \mathfrak{y}_k=1$ then $k\not\in S$ and $k\in T$ thence $S\ne T$.
\item If   $\mathfrak{x}_k=1 \wedge \mathfrak{y}_k=0$ then $k\in S$ and $k\not\in T$ thence $S\ne T$.
\end{itemize}

Hence
\[
\varphi(l_1)=(S,n)\ne(T,m)=\varphi(l_2).
\]
This means that $\varphi$ is a monomorphism.

(iv) Let $(S,n)\in L_2$ then define 
\[
\mathfrak{x}_i= \begin{cases} 1, &{\rm if}\; i\in S;\\
0, &{\rm if}\; i\not\in S. 
\end{cases}
\]
Now by definition of $\varphi$ we have $\varphi((\mathfrak{x}_i),n)=(S,n)$. Therefore $\varphi$ is surjective.

Thus $\varphi$ is an  isomorphism.
\rule{2mm}{2mm}

\begin{proposition}\label{prop10.10}
If $\mathcal{N}_2=\{((\mathfrak{x}_i),0)\,|\, (\mathfrak{x}_i)\in \underset{\mathbb{Z}}{\bigoplus} \; GF^+(2)\}$ and
\linebreak

$\mathcal{H}_2=\{(0,n)\,|\,n\in\mathbb{Z}\}$ then $\mathcal{N}_2\trianglelefteq L'_2$ and $\mathcal{H}_2\le L'_2$.
\end{proposition}

$\Box$ (i) Let $((\mathfrak{x}_i),0)\in\mathcal{N}_2$ then $((\mathfrak{x}_i),0)$ is the inverse element  of $((\mathfrak{x}_i),0)$. Let $((\mathfrak{y}_i),0)\in\mathcal{N}_2$ then
\[
((\mathfrak{x}_i),0)((\mathfrak{y}_i),0)=((\mathfrak{x}_i\oplus \mathfrak{y}_i),0)\in\mathcal{N}_2.
\]
Let $((\mathfrak{y}_i), m)\in \mathcal{N}_2((\mathfrak{x}_i),n)$ then $n=m$.
Let $((\mathfrak{x}_i),n)((\mathfrak{z}_i),0)=((\mathfrak{y}_i),n)$ then $\mathfrak{z}_i=\mathfrak{x}_i\oplus\mathfrak{y}_i$ and
\[
((\mathfrak{y}_i),n)=((\mathfrak{x}_i),n)((\mathfrak{x}_i\oplus\mathfrak{y}_i),0)\in((\mathfrak{x}_i),n)\mathcal{N}_2.
\]

(ii) Let $((0),n)\in\mathcal{H}_2$ then $((0),-n)$ is the inverse element  of $((0),n)$. Let  $((0),m)\in\mathcal{H}_2$ then
\[
((0),n)((0),m)=((0),n+m)\in\mathcal{H}_2. \quad
\rule{2mm}{2mm}
\]

\begin{proposition}
$L'_2$ is a semidirect product of $\mathcal{N}_2$ and $\mathcal{H}_2$.
\end{proposition}

$\Box$ Let $((\mathfrak{x}_i),n)\in L'_2$ then $((\mathfrak{x}_i),n)=((\mathfrak{y}_i),0)((0),n)$ where $\mathfrak{y}_i=\mathfrak{x}_{i+n}$. Thus $\mathcal{N}_2\mathcal{H}_2=L'_2$. Since $\mathcal{N}_2\cap\mathcal{H}_2=((0),0)$ then by \ref{prop10.10}. Proposition $L'_2$ is a  semidirect product of $\mathcal{N}_2$ and $\mathcal{H}_2$.
\rule{2mm}{2mm}

\begin{proposition}
If $\langle G, \cdot \rangle$ is a group then $\langle G, \circ \rangle$ is a group where 
$x\circ y= yx$ and the groups $\langle G, \cdot \rangle$, $\langle G, \circ \rangle$
are anti-isomorphic. 
\end{proposition}

$\Box$ (i) Associative law: $x\circ(y\circ z)=(y\circ z)x=(zy)x=z(yx)=(x\circ y)\circ z$.

(ii) The neutral element of $\langle G, \circ \rangle$ is the neutral element of $\langle G, \cdot \rangle$.
Similarly the inverse element $x^{-1}$ of $x$ in  $\langle G, \circ \rangle$ is the inverse element in $\langle G, \cdot \rangle$. 

(iii) The map $\mathbb{I}: x\mapsto x$ is a bijection,  furthermore
$\mathbb{I}(xy)=xy=y\circ x=\mathbb{I}(y)\circ \mathbb{I}(x)$. Thus the groups $\langle G, \cdot \rangle$, $\langle G, \circ \rangle$
are anti-isomorphic. 
\rule{2mm}{2mm}

\begin{corollary}
The grupoid $L_2^a=\langle L_2, \circ \rangle$ is the group where 
\linebreak
$l_1\circ l_2= l_2l_1$, furthermore the groups $L_2$ and $L_2^a$ are anti-isomorphic.
\end{corollary}

\begin{definition}
The group $L_2^a$ is called the lamplighter group.
\end{definition}

\begin{corollary}
The grupoid $L_{\bar 2}^a=\langle L_2, \breve\circ \rangle$ is the group where 
\linebreak
$l_1 \breve\circ l_2= l_2\bar \circ l_1$, furthermore 
\begin{itemize}
	\item[\textbullet] the groups $L_{\bar 2}^a$ and $L_{\bar 2}$ are anti-isomorphic;
	\item[\textbullet] $L_{\bar 2}^a\cong L_2^a$.
\end{itemize}
\end{corollary}

\begin{corollary}\label{c10.14}
The grupoid $L''_2=\langle L'_2, \circ \rangle$ is the group where 
\linebreak
$l_1\acute\circ l_2= l_2l_1$, furthermore 
\begin{itemize}
	\item[\textbullet] the groups $L''_2$ and $L'_2$ are anti-isomorphic;
	\item[\textbullet] $L''_2\cong L_2^a\cong L_{\bar 2}^a$.
\end{itemize}

\end{corollary}

\begin{lemma}
$\mu[f^n]\alpha[h]\mu[f^{-n}]=\alpha[f^{-n}h]$
\end{lemma}

$\Box \qquad \mu[f]\alpha[h]\mu[f^{-1}]\underset{\ref{prop7.7}.\rm{Prop}}{=}\alpha[f^{-1}h]$

The rest is induction.
\begin{eqnarray*}
\mu[f^{n+1}]\alpha[h]\mu[f^{-n-1}] &\underset{\ref{l8.4}.\rm{Lem}}{=}& \mu[f]\mu[f^n]\alpha[h]\mu[f^{-m}]\mu[f^{-1}]\\ \\
&=& \mu[f]\alpha[f^{-n}h]\mu[f^{-1}] \\ 
&\underset{\ref{prop7.7}.\rm{Prop}}{=}&\alpha[f^{-1}f^{-n}h]=\alpha[f^{-n-1}h] \quad \rule{2mm}{2mm}
\end{eqnarray*}

\begin{lemma}
If $\alpha_n=\alpha[f^{n_1}]\alpha[f^{n_2}]\ldots\alpha[f^{n_j}]$ then
\[
\mu[f^\varkappa]\alpha_n\mu[f^{-\varkappa}]=\alpha[f^{n_1-\varkappa}]\alpha[f^{n_2-\varkappa}]\ldots\alpha[f^{n_j-\varkappa}]
\]
\end{lemma}

$\Box$
\begin{eqnarray*}
 \mu[f^\varkappa]\alpha_n\mu[f^{-\varkappa}]
&=&\mu[f^\varkappa] \alpha[f^{n_1}]\alpha[f^{n_2}]\ldots\alpha[f^{n_j}] \mu[f^{-\varkappa}]\\
&=&\mu[f^\varkappa] \alpha[f^{n_1}] \mu[f^{-\varkappa}] \mu[f^\varkappa] \alpha[f^{n_2}] \mu[f^{-\varkappa}] \mu[f^\varkappa]\ldots \\
&\ldots&\mu[f^{-\varkappa}] \mu[f^\varkappa]\alpha[f^{n_j}] \mu[f^{-\varkappa}]\\
&=&\alpha[f^{-\varkappa}f^{n_1}]\alpha[f^{-\varkappa}f^{n_2}]\ldots\alpha[f^{-\varkappa}f^{n_j}]\\
&=&\alpha[f^{n_1-\varkappa}]\alpha[f^{n_2-\varkappa}]\ldots\alpha[f^{n_j-\varkappa}] \quad \rule{2mm}{2mm}
\end{eqnarray*}

\begin{theorem}
$\Gamma(\mathcal{M}_2)\cong L^a_{2}$
\end{theorem}

$\Box$ Let
\begin{eqnarray*}
\alpha_k &=&\alpha[f^{k_1}]\alpha[f^{k_2}]\ldots\alpha[f^{k_i}],\\
\alpha_n &=& \alpha[f^{n_1}]\alpha[f^{n_2}]\ldots\alpha[f^{n_j}]
\end{eqnarray*}
then \qquad $(\alpha_k,\mu[f^\varkappa])\oslash (\alpha_n,\mu[f^\nu])=$
\begin{eqnarray*}
 &=& (\alpha_k\mu[f^\varkappa]\alpha_n\mu[f^{-\varkappa}],\mu[f^\varkappa]\mu[f^\nu])\\
&=&(\alpha_k\alpha[f^{n_1-\varkappa}]\alpha[f^{n_2-\varkappa}]\ldots\alpha[f^{n_j-\varkappa}],\mu[f^{\varkappa+\nu}])
\end{eqnarray*}
and by \ref{teor9.3}. Theorem $\Gamma(\mathcal{M}_2)\cong\langle N\times H, \oslash\rangle$ where 
\[
N=\langle \alpha[sf^m]\,|\,s\in\{0,1\} \wedge m\in\mathbb{Z}\rangle, \quad \; H=\langle \mu[f]\rangle.
\]

Let $\phi:\langle N\times H, \oslash\rangle\to L^a_{\bar 2}:(\alpha_k,\mu[f^\varkappa])\mapsto(\phi_1(\alpha_k),\phi_2(\mu[f^\varkappa]))=(S,\varkappa)$ where $S=\{k_1,k_2,\ldots, k_i\}$. 

{\bf Note.}\quad  $\phi$ is a bijection and 
\begin{eqnarray*}
&&\phi((\alpha_k,\mu[f^\varkappa])\oslash(\alpha_n,\mu[f^\nu]))\\
&=& \phi(\alpha_k\alpha[f^{n_1-\varkappa}]\alpha[f^{n_2-\varkappa}]\ldots\alpha[f^{n_j-\varkappa}],\mu[f^{\varkappa+\nu}])\\
&=& (\phi_1(\alpha_k\alpha[f^{n_1-\varkappa}]\alpha[f^{n_2-\varkappa}]\ldots\alpha[f^{n_j-\varkappa}]),\phi_2(\mu[f^{\varkappa+\nu}]))\\
&=&(\phi_1(\alpha_k\alpha[f^{n_1-\varkappa}]\alpha[f^{n_2-\varkappa}]\ldots\alpha[f^{n_j-\varkappa}]),\varkappa+\nu)
\end{eqnarray*}
Otherwise
\begin{eqnarray*}
&&\phi(\alpha_k,\mu[f^\varkappa])\breve\circ\phi(\alpha_n,\mu[f^\nu])=(S,\varkappa)\breve\circ(T,\nu)=(S\bigtriangleup T_{-\varkappa},\varkappa+\nu),
\end{eqnarray*}
where $T=\{n_1,n_2,\ldots, n_j\}$.

We must prove that 
\[
\phi_1(\alpha_k\alpha[f^{n_1-\varkappa}]\alpha[f^{n_2-\varkappa}]\ldots\alpha[f^{n_j-\varkappa}])=S\bigtriangleup T_{-\varkappa}.
\]

{\bf Note.}\quad  $\phi_1(\alpha_k)=S$ and
\[
\phi_1(\alpha[f^{n_1-\varkappa}]\alpha[f^{n_2-\varkappa}]\ldots\alpha[f^{n_j-\varkappa}])=\{n_1-\varkappa,n_2-\varkappa,\ldots,n_j-\varkappa\}=T_{-\varkappa}.
\]

 If $\exists \tau\; \exists t\; \tau=n_t-\varkappa$ then $\alpha[f^\tau]\alpha[f^{n_t-\varkappa}]=0$.
We know $\forall\tau\; \forall t\;\alpha[f^\tau]\alpha[f^t]=\alpha[f^t]\alpha[f^\tau]$. Therefore
\[
\phi_1(\alpha_k\alpha[f^{n_1-\varkappa}]\alpha[f^{n_2-\varkappa}]\ldots\alpha[f^{n_j-\varkappa}])=(S\cup T_{-\varkappa})\setminus(S\cap T_{-\varkappa})=S\bigtriangleup T_{-\varkappa}.
\]
Hence $\langle N\times H, \oslash\rangle\cong L_{\bar 2}^a$. Thereby $\Gamma(\mathcal{M}_2)\cong\langle N\times H, \oslash\rangle \cong L_{\bar 2}^a\cong L_2^a$. \rule{2mm}{2mm}
\bigskip

So up to isomorphism $\Gamma(\mathcal{M}_2)$ is the lamplighter group.

\end{document}